\documentclass[invmat,twoside,10pt]{svjour}
\usepackage{amsmath,amssymb,amscd}
\usepackage{graphicx}
 \DeclareGraphicsExtensions{.eps}
\usepackage{mathrsfs}
\usepackage[mathcal]{eucal}

\spnewtheorem{Thm}{Theorem}{\bfseries}{\itshape}
\spnewtheorem{Cor}{Corollary}{\bfseries}{\itshape}
\spnewtheorem{Prop}[Cor]{Proposition}{\bfseries}{\itshape}
\spnewtheorem{Lem}[Cor]{Lemma}{\bfseries}{\itshape}
\spnewtheorem{PLem}[Cor]{Principal Lemma}{\bfseries}{\itshape}
\spnewtheorem{Conj}[Cor]{Conjecture}{\bfseries}{\itshape}
\spnewtheorem{Def}[Cor]{Definition}{\bfseries}{\rmfamily}
\spnewtheorem{Ex}[Cor]{Example}{\scshape}{\rmfamily}
\spnewtheorem{Rem}[Cor]{Remark}{\scshape}{\rmfamily}

\makeatletter
 \def\subsubsection{\@startsection{subsubsection}{3}{\z@}%
    {-13dd plus-4pt minus-4pt}{-5.5pt}{\itshape}}
\makeatother

\renewcommand\ge{\geqslant}
\renewcommand\le{\leqslant}
\let\tildeaccent=\~
\let\hataccent=\^
\renewcommand\~[1]{\widetilde{#1}}
\renewcommand\^[1]{\widehat{#1}}

\def\<{\left<}
\def\>{\right>}
\def\({\left(}
\def\){\right)}

\def\abs#1{\left\vert #1 \right\vert}
\def\size#1{\mathbf S\left(#1 \right)}

\let\parasymbol=\S
\def\secref#1{\parasymbol\ref{#1}}
\def\ssm{{\,\text{--}\,}}
\def\pd#1#2{\frac{\partial#1}{\partial#2}}
\let\log=\ln

\let\simeq=\cong

\let\subsetneq=\subsetneqq

\let\polishL=l
\def\Zoladek.{\.Zol\c adek}

\def\Var{\operatorname{Var}}
\def\const{\operatorname{const}}

\def\Mat{\operatorname{Mat}}

\def\Re{\operatorname{Re}}

\def\Arg{\operatorname{Arg}}
\def\dist{\operatorname{dist}}

\def\ord{\operatorname{ord}}
\def\GL{\operatorname{GL}}
\def\etc.{\emph{etc}.}

\def\:{\colon}
\def\R{{\mathbb R}}
\def\C{{\mathbb C}}
\def\Z{{\mathbb Z}}
\def\N{{\mathbb N}}
\def\Q{{\mathbb Q}}
\def\P{{\mathbb P}}
\def\I{{\mathbb I}}
\def\e{\varepsilon}
\def\S{\varSigma}
\def\l{\lambda}
\let\ol=\overline
\def\f{\varphi}
\def\F{{\varPhi}}
\def\iu{\mathrm{i}}
\def\diag{\operatorname{diag}}
\def\conf{\asymp}
\def\poly{\operatorname{\textup{Poly}}}

\def\d{\,\mathrm d}

\def\G{\varGamma}
\let\PolishL=\L
\def\Lojas.{\PolishL ojasiewicz}
\def\L{\varLambda}
\def\cN{{\mathcal N}}
\def\scS{{\mathscr S}}
\def\scC{{\mathscr C}}
\def\Aut{\operatorname{Aut}}
\def\cdiam#1{\operatorname{cdiam}(#1)}
\def\slope{\sphericalangle}

\def\bullet{{\,\boldsymbol\cdot}\,}
\def\lbar{{\hat\l}}

\let\stitem=\item
\def\item[#1]{\stitem[\textup{#1}]}
\def\manualqed{\relax}

\begin{document}

\title { On the number of zeros of Abelian integrals }
\subtitle{A constructive solution of the Infinitesimal Hilbert
Sixteenth Problem}

\author
 {Gal Binyamini \and
  Dmitry Novikov \and
  Sergei Yakovenko
 }

\dedication{To Yulij Sergeevich Ilyashenko, who discovered this
problem 40 years ago, for his 65th birthday with gratitude and
admiration.}

\institute{Weizmann Institute of Science, Rehovot, Israel\\
\email{\{gal.binyamini,\,dmitry.novikov,\,sergei.yakovenko\}@weizmann.ac.il}}

\mail{Sergei Yakovenko}

\date{Posted on August 21, 2008. Revised version Oct 5, 2009.}

\maketitle

\begin{abstract}
We prove that the number of limit cycles generated from nonsingular
energy level ovals (periodic trajectories) in a small
non-conservative perturbation of a Hamiltonian polynomial vector
field on the plane, is bounded by a double exponential of the degree
of the fields. This solves the long-standing infinitesimal Hilbert
16th problem.

The proof uses only the fact that Abelian integrals of a given
degree are horizontal sections of a regular flat meromorphic
connection  defined over $\Q$ (the Gauss-Manin connection) with a
quasiunipotent monodromy group.
\end{abstract}

\keywords{Abelian integrals -- Fuchsian systems -- monodromy --
limit cycles. \\[5pt]
}

\setcounter{tocdepth}{2} \tableofcontents

\setitemindent{~(Q)~}


\section{Infinitesimal Hilbert 16th problem}

The central result of this paper is an explicit upper bound for the
number of limit cycles born from nonsingular (smooth compact) energy
level ovals in a non-conservative polynomial perturbation of a
polynomial Hamiltonian vector field on the plane. This problem was
repeatedly posed in various sources under different names as the
weakened, infinitesimal or tangential Hilbert 16th problem. In this
introductory section we briefly outline some connections between
different problems concerning limit cycles of polynomial vector
fields on the plane. Much more complete expositions can be found in
the recent surveys \cite{centennial,montreal} and the books
\cite{Zol,iy:lade}.

\subsection{Limit cycles born by perturbations of
 integrable planar vector fields}\label{sec:perturb-cycles}
Limit cycles, \emph{isolated} periodic (compact, nontrivial)
trajectories of polynomial vector fields, are one of the most
elusive objects of analysis. There are only a handful of tools to
establish the (non)existence of such cycles in certain domains, all of
them applying to very specific differential equations. D.~Hilbert,
who included the problem of counting the possible number of limit
cycles as the 16th problem in his famous list \cite{hilbert},
conjectured implicitly that the problem could be approached by
perturbation techniques, first studying vector fields close to those
which are ``simple'' from the point of view of counting their limit
cycles. One such natural class is the \emph{integrable fields} which
exhibit continuous families of (nonisolated) periodic trajectories
(here and below \emph{integrability} means existence of a local or
global \emph{first integral}). In  \cite{petr-land:deg=n}
Petrovski{\u\i} and Landis attempted to realize this program by
a complexification of the problem, but their attempt was not successful
\cite{ilyash:mathevents}. Nevertheless, the problem of estimating
the number of limit cycles of near-integrable systems became a
natural intermediate step towards a possible future solution of the
Hilbert problem which still seems to be completely out of reach,
see \cite[Problem 1978-6, pp.~352--363]{arnold-problems} and, most
recently, \cite[\parasymbol 3.2]{ilyash:openprob}.

Among the integrable systems the simplest (in many respects) class
is that of \emph{Hamiltonian polynomial systems}, vector fields
corresponding to a system of autonomous differential equations of
the form
\begin{equation}\label{hamsys}
    \frac{\d x}{\d z}=\pd Hy (x,y),\qquad \frac{\d y}{\d z}=-\pd Hx
    (x,y),
\end{equation}
with a real bivariate polynomial $H$ called the Hamiltonian; in
the Pfaffian form these differential equations can be written
as $\d H=0$, and all real level ovals of $H$ (compact connected
components of the level curves of the form $\{H=t\}$) are
integral trajectories.

A polynomial perturbation of the Hamiltonian system
\eqref{hamsys} can also be written in the Pfaffian form with a
small parameter $\e\in(\R^1,0)$ as follows,
\begin{equation}\label{ham-pert}
    \d H+\e\omega=0,\qquad  \omega=P(x,y)\d x+Q(x,y)\d y.
\end{equation}
In general, such perturbations destroy
integrability, so that for $\e\ne0$ most integral trajectories will
become spirals.

We say that a (smooth) closed oval $\delta\subseteq\{H=t\}$
generates a limit cycle in the perturbation \eqref{ham-pert}, if for
any sufficiently small annular neighborhood $U$ of $\delta$ one can
find arbitrarily small values of the parameter $\e$ such that the
corresponding Pfaffian equation exhibits a limit cycle $\delta_\e$
entirely belonging to $U$. If there exists a natural number $k\ge 1$
such that for an arbitrarily narrow $U$ and an arbitrarily small
$\e$ there may coexist $k$ limit cycles, we say that the oval
$\delta$ generates $\ge k$ limit cycles in the family.

Respectively, we say that an oval $\delta$ generates no more
than $k$ limit cycles in perturbation \eqref{ham-pert}, $k\ge
0$, if there exists a small annular neighborhood $U$ of
$\delta$ on the $(x,y)$-plane, and a small neighborhood
$V=(\R^1,0)$ of the origin on the parameter axis, such that for
any $\e\in V$ the foliation defined by the Pfaffian form $\d
H+\e\omega$, has no more than $k$ limit cycles in $U$. The
minimal number $k$ with this property, denoted by
$k=k(\delta;\omega)$, always exists: one can easily see that
for almost all ovals $k(\delta;\omega)=0$ (the oval is
destroyed without generating any limit cycle). In other words,
in the sum taken over all smooth ovals of the real level curves
$\{H=\const\}\subset\R^2$,
\begin{equation}\label{total-num}
    \cN(H,\omega)=\sum_{\delta\subseteq\{H=\const\}}k(\delta;\omega)\le+\infty
\end{equation}
all but countably many terms mush vanish.

It is well known that $k(\delta;\omega)>0$ only if the
Poincar\'e integral
\begin{equation}\label{poincare-int}
    I=I(\delta,\omega)=\oint_{\delta}P\d x+Q\d y,
\end{equation}
vanishes (Poincar\'e--Andronov--Pontryagin criterion
\cite[\parasymbol 26A]{iy:lade}\footnote{For convenience of the
reader we give references to the textbook \cite{iy:lade} whenever
possible. References to the original publications can be found in
this textbook.}). In physical terms, the integral
\eqref{poincare-int} is the principal asymptotic term for the
dissipation of the energy along one period. The perturbation is
called \emph{non-conservative} if $I(\cdot,\omega)\not\equiv0$. In
this case a slight refinement of the
Poincar\'e--Andronov--Pontryagin criterion asserts that
$k(\delta;\omega)$ does not exceed the multiplicity of the root of
the integral $I(\delta,\omega)$ as the function of the first
argument.

\medskip

The infinitesimal Hilbert problem requires to \emph{place an upper
bound for the number of limit cycles born from nonsingular energy
level ovals in non-conservative perturbations}. The answer should
depend only on $n$, i.e., the bound must be uniform over all
Hamiltonians $H$ of degree\footnote{The degrees of the Hamiltonian
and the perturbation $1$-form are chosen so that both terms in the
perturbation \eqref{ham-pert} have the same degree $n$.} $\deg H\le
n+1$ and polynomial $1$-forms $\omega$ of degree $\deg
\omega=\max(\deg P,\deg Q)\le n$. Our main result solves this
problem and gives an explicit double exponential upper bound.

\begin{Thm}\label{thm:main}
The total number of limit cycles $\cN(H,\omega)$ that can be
born from {all nonsingular} energy level ovals of a Hamiltonian
polynomial foliation in a non-conservative perturbation
\eqref{ham-pert} of degree $\le n$, is no greater than
$2^{2^{\poly(n)}}$.

Here the expression $\poly(n)=O(n^p)$ stands for an explicit
polynomially growing term with the exponent $p$ not exceeding $61$.
\end{Thm}

Besides limit cycles born from nonsingular ovals, limit cycles
can be born from separatrix polygons (energy level curves
carrying singular points of the Hamiltonian vector field).
Theorem~\ref{thm:main} does not address the number of these
cycles, see Remark~\ref{rem:separatrix} below.

\subsection{Zeros of Abelian integrals}
%
%
The Poincar\'e integral for polynomial perturbations is an
integral of a rational (in fact, polynomial) 1-form $\omega$
over a cycle $\delta$ on the real algebraic curve $\{H=t\}$.
Such integrals are called \emph{Abelian integrals} and they can
be considered as functions of all parameters occurring in the
construction (coefficients of the 1-form and the algebraic
curve). In particular, we can consider the Hamiltonian $H$ and
the 1-form $\omega$ as the parameters and look at the
Poincar\'e integral \eqref{poincare-int} as a continuous branch
of a multivalued function $I_{H,\omega}(t)$ of the single
variable $t$, the value taken by $H$ on the cycle
$\delta(t)\subseteq\{H(x,y)=t\}$\footnote{If the curve
$\{H=t\}$ carries several smooth real ovals, we always consider
several branches of the integral $I$ simultaneously.}. For a
non-conservative perturbation, this function is not identically
zero by definition.

In such a context the question about limit cycles born in the
perturbation \eqref{ham-pert} reduces to the question about  an
upper bound for the \emph{total number of isolated zeros of
real branches of the Abelian integral $I_{H,\omega}(t)$,
counted with multiplicities.} This upper bound should be
uniform over arbitrary combinations of the ``parameters'' $H$
and $\omega$ of degrees not exceeding a given natural number
$n$.

Theorem~\ref{thm:main} is an immediate corollary of the
following result on zeros of Abelian integrals. For any finite
values $n,m\in\N$ denote by $\cN(n,m)$ the upper bound for the
number of isolated real zeros of the integrals $I$ counted with
multiplicities,
\begin{multline}\label{vk-fin}
    \cN(n,m)=\sup_{\omega,H}\cN(H,\omega)=
    \sup_{\omega,H}~\sum_t\ord_t I_{H,\omega}(\bullet),\\ \deg H\le n+1,\
    \deg\omega\le m.
\end{multline}
Here $\ord_t I\ge0$ denotes the order of the root of the integral
$I$ at a real point $t$ (this order is zero by definition if
$I(t)\ne 0$ and 1 for a simple root); if for a given value of $t$
the algebraic curve $\{H=t\}$ carries several real nonsingular
ovals, the summation is extended over all corresponding continuous
branches of $I$.
%

\begin{Thm}\label{thm:main2} $$\cN(n,n)\le 2^{2^{\poly(n)}},$$
were the expression $\poly(n)=O(n^p)$ stands for an explicit
polynomially growing term with the exponent $p$ not exceeding $61$.
\end{Thm}

This result, the first \emph{explicit uniform bound} for the number
of isolated zeros of Abelian integrals, is the most recent
(hopefully, not the last) step in the long line of research, partly
outlined in \secref{sec:milestones}.

\begin{Rem}[notation for polynomial bounds]
 \label{rem:oplus}
Here and below we will deal with explicit bounds which involve
simple or double exponentials of polynomially growing terms. In
order to avoid cumbersome notation, we will use the following
shortcuts.

Everywhere unless explicitly waived, the symbol $O(\bullet)$ refers
to an explicit constructive asymptotic; in particular, the notation
$O(n^p)$ means a constructive positive function of an integer
argument $n$ which does not exceed the expression $Cn^p$ for some
explicit constant $C<+\infty$ and all $n\ge 2$.

The notation $\poly(n)$ stands for the constructive bound $O(n^p)$
for some unspecified finite exponent $p<+\infty$.

Sometimes the growth rate will be estimated by the ``extended
polynomial notation'' $O^+(n^p)$, which by definition means
``$O(n^{p+\e})$ for any positive $\e$''. A typical example is
as follows: for any finite $q>0$,
\begin{equation}
    n\log^q n=O^+(n).
\end{equation}
Of course, the notation $O^+(n^p)$ in the upper bound could be
replaced by $O(n^{p+1})$, yet when such ``rounding errors'' are
composed, the overall bound gets increased by artificial terms
unrelated to the nature of these bounds.

For bounds depending on several arguments, we use the notation
$O(n^pm^q)$, $\poly(n,m)$ and $O^+(n^pm^q)$ in a self-explanatory
way.

The introduced notation allows to formulate the improved bounds in
Theorems~\ref{thm:main} and~\ref{thm:main2}. What we prove in fact
is the double exponential bound $2^{2^{\poly(n)}}$ with
$\poly(n)=O^+(n^{60})$, see \eqref{oplus-main}.
\end{Rem}

\begin{Rem}[accuracy of the upper bounds]
The double exponential expression cited in Theorem~\ref{thm:main2}
is only an \emph{upper bound} with no claim of accuracy whatsoever.
Moreover, some fine tuning of the tools developed in this paper, can
apparently help in reducing the power $p=61$ to a much lower value,
perhaps, as low as $p=2$ (at the price of clarity of the
exposition), see \secref{sec:concl}. Reducing the bounds to less
than two exponential orders would definitely require new ideas.

On the other hand, the only known \emph{lower bounds} are
quadratic in $n$ and linear in the degree of the form
$m=\deg\omega$ for a fixed $n$ if the latter is allowed to grow
to infinity independently of the degree of the Hamiltonian,
cf.~with \eqref{vk-fin}. Thus far there is no sound conjecture
on what might be the true asymptotic behavior of the function
$\cN(n,m)$.
\end{Rem}

\subsection{Bifurcations of limit cycles \emph{not covered} by
Theorem~\ref{thm:main}}\label{sec:bifurcations}
 \begin{small}
In order to avoid possible ambiguities, we  make several remarks on
what \emph{is not asserted} in Theorem~\ref{thm:main}. The remarks
below can be considered as a list of open problems.\par

\begin{Rem}[Conservative and integrable perturbations]
Theorem~\ref{thm:main} gives no bound on the number of limit cycles if
the perturbation itself is \emph{conservative}, i.e., if the Poincar\'e
integral vanishes \emph{identically} for all ovals
$\delta$ in the family. For instance, if $\omega$ is exact,
$\omega=\d F$, $F\in\R[x,y]$, then the entire family
\eqref{ham-pert} consists of Hamiltonian systems and exhibits no
limit cycles for all $\e\in(\R,0)$.

On the other hand, the identical vanishing of the Poin\-car\'e
integral \eqref{poincare-int} in general \emph{does not imply}
that all foliations in the family \eqref{ham-pert} are
integrable and hence do not have limit cycles at all. Indeed,
the integral \eqref{poincare-int} is only the first variation
of the Poincar\'e return map with respect to the parameter
$\e$. If the first variation vanishes, one can compute higher
variations in what is sometimes called the \emph{Fran\c coise
algorithm} \cite{francoise,bautin}, see also \cite[\parasymbol
26B]{iy:lade}. The number of isolated zeros of the \emph{first
not identically vanishing variation} will majorize the number
of limit cycles born from smooth ovals of the Hamiltonian
field.

For a \emph{generic} polynomial $H$, one can show that the higher
variations will again be Abelian integrals of certain polynomial
1-forms, yet their degrees are growing together with the order of
the variation. Thus the problem of counting the limit cycles in the
perturbation \eqref{ham-pert} reduces via Theorem~\ref{thm:main} to
the question on \emph{how many consecutive higher variations can
vanish identically without the family \eqref{ham-pert} being
necessarily integrable}. This is a generalization of the famous
\emph{Poincar\'e center problem} whose solution is unknown even in
the best of the best cases $H(x,y)=x^2+y^2$.

In the special degenerate (e.g., symmetric) cases the higher
variations cannot be expressed as Abelian integrals, only
by means of (polynomial expressions in) the so called \emph{iterated integrals}
\cite{gavrilov-iliev:iterated,gavrilov:iterated}. Formally these
integrals are not covered by Theorem~\ref{thm:main2} below. However,
L.~Gavrilov and I.~Iliev have shown that the iterated integrals still
satisfy a Fuchsian system of equations and very recently it was
discovered that the monodromy group of this system is quasiunipotent
\cite{novikov-gavrilov-2008}. These observations pave a way to
application of Theorem~\ref{thm:count}, our principal result, yet
many things remain to be done in order to bridge the gaps. Anyhow, the
same question on the maximal order of the nontrivial high variations
(an equivalent of the center problem) reappears in this case as
well.
\end{Rem}

\begin{Rem}[Various scenarios of integrability]
Hamiltonian vector fields are only the simplest case of integrable
polynomial foliations on the projective plane $\R P^2$. The question
about all possible scenarios of integrability is one notch above the
Poincar\'e center problem (which addresses the question of
integrability only locally, near a singular point), hence is
challenging and wide open. Still, some possibilities are well known
and documented.

One such scenario is the \emph{Darbouxian integrability}, which
corresponds to replacing the \emph{exact polynomial} 1-form $\d H$
in \eqref{ham-pert}, by a \emph{closed rational} 1-form (note that
the equations in the Pfaffian form admit multiplication by a
rational factor without changing the behavior of the trajectories).

Investigation of limit cycles born by perturbation of Darbouxian
integrable systems is a completely new field, where only the first
steps are now taken \cite{mit:arxiv:Darboux,bob-mar-JLond}. One of
the main difficulties is the fact that the ovals of Darbouxian
integrable systems are in general not algebraic. This circumstance
renders practically inapplicable all known tools working for Abelian
integrals. In particular, Theorem~\ref{thm:main} seems to be of no
help in this context, as the Poincar\'e integrals do not satisfy
any finite order linear differential equation.

Besides Hamiltonian and Darbouxian integrability, there are some
other known scenarios. The most difficult for analysis seems to be
appearance of limit cycle by perturbation of \emph{symmetric
systems}, yet the problem is too vague to be discussed here.
\end{Rem}

\begin{Rem}[Limit cycles born from nonsmooth level curves of the
 Ham\-ilt\-onian]\label{rem:separatrix}
Theorem~\ref{thm:main} asserts \emph{nothing} about the number of limit
cycles born from \emph{nonsmooth} ovals, corresponding to the
\emph{critical level curves} of Hamiltonians.

For a \emph{generic} real Hamiltonian, the only critical level ovals
are separatrix loops (homoclinic trajectories of a nondegenerate
saddle) and double loops (eight-shaped curves, \emph{butterflies}),
pairs of homoclinic orbits of the same saddle, which may generate
limit cycles converging uniformly to the union of the two
trajectories.

The case of separatrix loops is well understood after the seminal
works of R.~Roussarie \cite{roussarie:brasil,roussarie:cyclicite}.
An upper bound for the number of cycles born from a simple loop can
be obtained from Theorem~\ref{thm:main2}, yet much weaker results on
the maximal multiplicity of zero of Abelian integrals, in the spirit
of \cite{mardesic,moura:jdcs} are sufficient and give much better
bounds. A similar bound for the double loops could perhaps be
derived using the tools from \cite{jebrane-mourtada}.

For singular level curves carrying \emph{more than one} saddle
critical point of $H$, one cannot in general predict the number of
limit cycles based only on the first asymptotic terms of the Abelian
integral (a substitute for the order of zero for points of
non-analyticity of $I$). In \cite{dum-rous:alien,cdr:alien} it is
shown that already for perturbations of a Hamiltonian foliation with
two nondegenerate critical points on the same level curve, one can
obtain limit cycles not related to zeros of Abelian integrals
(called \emph{alien cycles}).

The general question about limit cycles born from a critical
level curve of an arbitrary polynomial Hamiltonian, is quite
challenging and essentially open.
\end{Rem}
\par
\end{small}

 \subsection{A few milestones}\label{sec:milestones}
 \begin{small}
Probably the first to realize that investigation of limit
cycles in near-integrable systems may be a path to solution of
the Hilbert problem,  were I. G. Petrovski\u\i\ and E. M.
Landis \cite{petr-land:deg=2,petr-land:deg=n,petr-land:corr}.
Although their direct approach turned out to be unfeasible,
these seminal papers stimulated the study of perturbations of
Hamiltonian systems.

The first study of Abelian integrals, focused on investigation of
their roots in connection with the bifurcation of limit cycles, was
undertaken in the dissertation of Yu.~Ilyashenko (adviser E. M.
Landis), see \cite{ilyashenko:1969,ilyashenko:1969b}. In this work
Ilyashenko introduced very powerful tools from complex analysis and
algebraic geometry and implicitly formulated the infinitesimal
Hilbert problem in the form we solve it now.

One of the first explicit bounds for the number of zeros of Abelian
(elliptic) integrals was obtained by R.~Bogdanov \cite{bogdanov:ai};
soon Yu.~Ilyashenko gave a transparent proof this result
\cite{ilyashenko:elliptic}.

Since then the number of papers devoted to investigation of zeros of
Abelian integrals counts in the hundreds, and it is impossible to
mention even the names of the principal contributors. The
overwhelming majority of these papers deal with the low-degree cases
$n=3,4$, where the problem is essentially settled by L.~Gavrilov,
I.~Iliev and C.~Li (see part 2 of the book \cite{li-christopher} and
references therein).

In the general case of arbitrary degree, however, very little is
known. In 1984 A.~Khovanski{\u\i} and A.~Varchenko proved that the
number of zeros of Abelian integrals is always finite and uniformly
bounded over all Hamiltonians and 1-forms of a given degree
\cite{asik:finiteness,varchenko:finiteness}, i.e., that
$\cN(n,m)<+\infty$ for all finite combinations $n,m$. Unfortunately,
the proof is purely existential and does not give explicit bounds on
$\cN(n,m)$ even for small degrees $n,m$.

After this celebrated result many efforts were focused on obtaining
asymptotic constructive bounds for the ``counting function''
$\cN(\bullet)$. Very soon it became clear that the roles of the form
and the Hamiltonian are quite different from the point of view of
the difficulty of study. More precisely, for a fixed Hamiltonian $H$
one may consider forms of arbitrary degree $m$ growing to infinity,
and study the asymptotic of the counting function $\cN(H,m)$ as
$m\to+\infty$. The first bounds for this restricted setting were
double exponential in the degree of the form \cite{invmath-95}, yet
very quickly they were replaced by single exponential
\cite{annalif-95} and finally linear \cite{petrov:quartic-nonosc} bounds.
The ultimate result, due to Petrov and Khovanski\u\i, is the
following estimate,
\begin{equation}\label{zoladek-lin-bound}
    \forall n,m\in\N\qquad\cN(n,m)=O_n^{\text{exist}}(1)+O(m).
\end{equation}
Here the first term $O_n^{\text{exist}}(1)$ is a \emph{purely
existential} ``constant'' depending on $n$ (uniformly over all
Hamiltonians of degree $\le n+1$) and the second term $O(m)$
is, as usual, explicit and constructive. For quite some time
the proof of this result existed only in the oral tradition,
until it was published in the book \cite[Theorem 6.26]{Zol}.
The proof is based on the fact that the Abelian integrals of a
1-form of arbitrarily high degree $m$ over level ovals of a
Hamiltonian $H$ of degree $\le n+1$ can always be represented
as a linear combination of integrals of 1-forms of degree not
exceeding $2n$ with coefficients polynomial in $t$ of degree
$O(m/n)$ (see Theorem~\ref{thm:petrov-mod} below). Based on
this observation, one can \emph{conjecture} that there exists
an explicit (constructive) bound of the form
\begin{equation}\label{general-bound-conj}
    \cN(n,m)\le 2^{2^{\poly(n)}}+O(m)\qquad\text{as }n,m\to+\infty
\end{equation}
(Theorem~\ref{thm:main2} addresses only the bound for $\cN(n,n)$).
The proof could hopefully be obtained by a combination of the two
techniques, since our methods allow placing an explicit upper bound
on the first (existential) term in \eqref{zoladek-lin-bound}.

The most recent development in connection with the infinitesimal
Hilbert 16th problem is an explicit upper bound for the number of zeros
of Abelian integrals, uniform over all 1-forms of degree $\le n$,
finite for \emph{almost all} Hamiltonians $H$ of degree $n+1$, yet
\emph{non-uniform} in $H$. In a series of papers
\cite{glu:topology,glu:determinant,iglu:dan,iglu:restricted} A.
Glutsyuk and Yu.~Ilyashenko established this type of bound,
which grows exponentially as $H$ approaches an exceptional set of
Hamiltonians. This result was improved in \cite{polyfuchs} where a
bound growing \emph{polynomially} near the same exceptional set was
given.

The only class of Hamiltonians of arbitrarily high degree for which
uniform explicit upper bounds were previously known, is the class of
\emph{hyperelliptic Hamiltonians} of the form
$H(x,y)=\frac12y^2+Q(x)$, $Q\in\R[x]$. In \cite{era-99} it was
proved that the number of isolated zeros of hyperelliptic integrals
can be majorized by a \emph{tower function} (iterated exponent) of
$n=\deg Q$ under the technical assumption that all critical values
of the hyperelliptic Hamiltonian are real. However, the height of
this tower was much larger than 2 (corresponding to the double
exponent).

\par
\end{small}

\subsection{Hyperelliptic case}
The tools developed in this paper are sufficiently flexible to give
better results for some more specific classes of Abelian integrals.
For instance, if instead of a general bivariate polynomial $H(x,y)$
we consider only the \emph{hyperelliptic} Hamiltonians of the form
\begin{equation}\label{hyperel}
    H(x,y)=\frac12y^2+x^{n+1}+\l_1x^{n-1}+\cdots+\l_n,\qquad
    \l\in\R^n,
\end{equation}
then some steps of the proof can be skipped or improved, see
\secref{sec:concl}. As a result, we have a better bound for zeros of
the corresponding hyperelliptic integrals.

\begin{Thm}\label{thm:hyperel}
The number of complex isolated zeros of a hyperelliptic
Abel\-ian integral associated with the Hamiltonian
\eqref{hyperel}, is bounded by the explicit double exponent
$2^{2^{O^+(n)}}$.
\end{Thm}

This result, completely superseding the main result of our
previous work \cite{era-99}, is explained in
\secref{sec:concl}.

\begin{acknowledgement}
We are grateful to S. Basu, J. Bernstein, P. Deligne, V. Hinich and
N. Vorobjov for useful conversations, professional consultations and
illuminating remarks which turned out to be of key importance for
the entire construction.

Besides, we feel the need to acknowledge the pioneering work of our
colleagues and former collaborators, that is essentially if not
always directly used in the proof: G. Petrov, M. Roitman, A.
Grigoriev, V.~Arnold, A. Khovanski\u\i, A. Gabrielov, L.~Gavrilov,
H.~\Zoladek. and especially Yu. Ilyashenko, who stood at the origin
of the problem and suggested its particular case as a subject for
the M.Sc.~thesis research to one of the authors \cite{sy:thesis}.

Finally, we owe many thanks to the two referees for their most
attentive reading of the draft: the numerous remarks they made were
very helpful in bringing the exposition to its final form.

During the work on this project, two of the authors (G.B. and S.Y.)
were partially supported by the Minerva Foundation.  D.N.~gratefully
acknowledges the support of Soref New Scientists Start up Fund and
Fusfeld Research Fund. S.Y.~is the incumbent of the Gershon Kekst
Professorial Chair.
\end{acknowledgement}

\section{Background, settings, main constructions and strategy}

We begin by describing the (standard) complexification of the
Abelian integrals and reduce the infinitesimal Hilbert 16th problem
to a question about zeros of solutions to an integrable Pfaffian
system subject to a condition on its monodromy. The exposition in
this section primarily settles the context and notations for the
main body; the recent textbook \cite[\parasymbol 26]{iy:lade}
contains all necessary details.

\subsection{Complete Abelian integrals depending on parameters}
Let $\G=\G_{\l}\subset\P^2$ be the complex projective curve defined
in the affine chart $(x_1,x_2)$ on $\C^2\subset\P^2$ by the equation
\begin{equation}\label{G}
   H(x_1,x_2)=0, \qquad H=H(x,\l)=\sum_{0\le|\alpha|\le n+1}\l_\alpha x^\alpha
\end{equation}
(here and below we use the standard multiindex notation denoting by
$\alpha\in\Z_+^2$ the multiindex, $|\alpha|=\alpha_1+\alpha_2$,
$x^\alpha=x_1^{\alpha_1}x_2^{\alpha_2}$). The parameters $\l$ of
this equation naturally vary over the \emph{complex projective}
space $\P^m$ of dimension $m=\tfrac12(n+3)(n+2)-1=O(n^2)$.

For a generic combination of the parameters $\l$ the curve $\G_{\l}$
is a nonsingular (smooth) Riemann surface of genus $g=\frac12n(n-1)$,
transversal to the infinite line $\I=\P^2\ssm\C^2$. The (first)
homology group of $\G_\l\ssm\I$ in this case has the rank
$\ell=n^2$, see \cite[Theorem 26.31]{iy:lade}. The combination of
the parameters corresponding to the exceptional (non-smooth or
non-transversal to $\I$) curves $\G_\l$, is a projective algebraic
subset that will be denoted by $\S_*$. For an arbitrary point
$\l_*\notin\S_*$ one can fix a system of cycles
$\delta_1,\dots,\delta_\ell$ generating the homology
$H_1(\G_{\l_*},\Z)$ with integer coefficients and transport them
horizontally in the sense of the Gauss--Manin connexion. The result
is a multivalued framing of the fibers $H_1(\G_\l,\Z)$ associated
with the topological bundle $\G_\l\mapsto\l$ over $\P^m\ssm\S_*$,
ramified over $\S_*$ \cite[Corollary 26.28]{iy:lade}.

The cohomology of a generic fiber (curve) $\G_\l$ is generated by
restrictions of polynomial 1-forms on this curve. Let
$\omega_\alpha=\frac{x_1}{\alpha_1+1} \cdot x^\alpha \d x_2$ be
monomial  1-forms which are primitives of the 2-forms
$\mu_\alpha=x^\alpha\d x_1\land\!\d x_2$ with
$0\le\alpha_1,\alpha_2\le n-1$, i.e. $\d\omega_\alpha=\mu_\alpha$.
The number of such forms is exactly equal to $\ell=n^2$, and it is
known (see Appendix~\ref{sec:compl-pf} for details and references)
that the $\omega_\alpha$ generate the {cohomology} of a
\emph{generic} fiber $\G_\l$ with $\l\notin\S_*$ over $\C$. However,
for some exceptional fibers the forms $\omega_\alpha$ become linear
dependent after restriction on $\G_\l$: the corresponding set is a
proper algebraic subvariety, whose union with $\S_*$ will be denoted
by $\S\subset\P^m$.

\begin{Def}\label{def:period-X}
The \emph{period} matrix $X(\l)$ is the $(\ell\times\ell)$-square
(multivalued) analytic matrix function on $\P^m$,
\begin{equation}\label{X}
    X(\bullet)=\begin{pmatrix}
    \displaystyle\oint_{\delta_1(\bullet)}
    \kern-15pt\omega_1&\cdots&\displaystyle\oint_{\delta_\ell(\bullet)}\kern-15pt\omega_1
    \\
    \vdots&\ddots&\vdots
    \\
    \displaystyle\oint_{\delta_1(\bullet)}
    \kern-15pt\omega_\ell&\cdots&\displaystyle\oint_{\delta_\ell(\bullet)}\kern-15pt\omega_\ell
    \end{pmatrix},
\end{equation}
ramified over the locus $\S_*$ and nondegenerate  on $\P^m\ssm\S$.
\end{Def}

\begin{small}
\begin{Rem}[From projective spaces to pencils of
 lines]\label{rem:pencil}
In the formulation of Theorem~\ref{thm:main2} the Abelian integral
occurs as a function of a distinguished variable $t$, whereas all
other coefficients of the Hamiltonian $H$ are treated as parameters.
In the definition of the period matrix all coefficients of $H$ play
the same role. However, this difference is only superficial, and one
can consider $X(\l)$ as a parametric family of functions of one
(complex) variable, so that their isolated zeros can be counted.

If in the expression for $H$ in \eqref{G} all parameters $\lbar
=\{\l_\alpha:\alpha\ne(0,0)\}$ are fixed except for the free term
$\l_{0,0}$, we obtain a (complex projective) line
$L_{\lbar}\simeq\P^1$. These lines corresponding to different values
of $\lbar\in\P^{m-1}$ form a \emph{pencil of lines}, a family of
lines in $\P^m$ passing through the point in $\P^m$ with the
homogeneous coordinates $(0,\dots,0,1)$.

The space $\P^m$ is birationally equivalent to the product
$\P^{m-1}\times\P^1$. Moreover, the equivalence can be chosen
to map the lines from the above pencil to the lines of the form
$\P^1_{\lbar}=\{\lbar\}\times\P^1$. In suitable affine charts,
the equivalence is represented by the natural identification
$\C^{m}\simeq\C^{m-1}\times\C^1$, $\l\simeq(\lbar,t)$. After
such identification the period matrix \eqref{X} can be indeed
considered as a multivalued matrix function $X(\lbar,t)$
\emph{with singularities}, defined on $\P^1$ and depending on
the parameters $\lbar\in\P^{m-1}$ which vary over a compact
parameter space.\par
\end{Rem}
\end{small}

\subsection{Integrable Pfaffian systems with singularities}
Consider a smooth (nonsingular) algebraic variety $\mathbb M$  and
let $\Omega$ be a rational $(\ell\times\ell)$-matrix-valued 1-form
on $\mathbb M$ with a singular locus $\S\subset \mathbb M$ which is
an algebraic hypersurface (eventually, itself singular and
reducible).

The form $\Omega$ is \emph{integrable} if
$\d\Omega-\Omega\land\Omega=0$. The integrability condition is
necessary and sufficient for the \emph{local} existence of a
holomorphic nondegenerate matrix solution $X(\bullet)$ for the
Pfaffian system of equations
\begin{equation}\label{pfs}
    \d X=\Omega \cdot X,\qquad\Omega\in\Mat\bigl(\ell,\L^1(\mathbb M)\bigr)
\end{equation}
near each nonsingular point $a\notin\S$. The local solution admits
analytic continuation along any path $\gamma$ avoiding the singular
locus. If the path $\gamma$ is closed, the result of continuation
$\Delta_\gamma X$ differs from the initial solution $X$ by a
constant nondegenerate matrix factor $M_\gamma\in\GL(\ell,\C)$,
called the \emph{monodromy} associated with this path:
$\Delta_\gamma X=XM_\gamma$.

\begin{Def}\label{def:small-loop}
Let $\tau\:(\C^1,0)\to(\mathbb M,\S)$, be the germ of a
holomorphic curve, not entirely belonging to $\S$, and
$a=\tau(0)$.

A \emph{small loop around the point $a$} is a closed path which is
the image $\tau(S_\e)$ of any sufficiently small circle
$S_\e=\{|s|=\e\}\subset(\C,0)$. Here the smallness means that the
image of the punctured disk $\{0<|s|\le \e\}$ is disjoint with $\S$.
\end{Def}

A loop freely homotopic to a small loop will also be referred to as
the small loop. All small loops ``supported'' by the same
holomorphic curve $\tau$ are freely homotopic to each other. The
integrability assumption implies that the corresponding monodromy
operators are conjugate to each other, and in particular they have
the same spectrum.

\begin{Def}\label{def:quasi-unipotent}
The integrable form $\Omega$ is called \emph{quasiunipotent} at a
point $a\in\mathbb M$, if all eigenvalues of each monodromy operator
associated with any small loop around $a$, are roots of unity.
Clearly, \emph{any} integrable form is quasiunipotent at a
nonsingular point $a\notin\S$, since the corresponding monodromy is
the identity. The form is (globally) \emph{quasiunipotent}, if it
is quasiunipotent at every  point of the singular locus of $\Omega$.
\end{Def}

\begin{Rem}\label{rem:small-loops-1d}
If $\dim\mathbb M=1$, i.e., the system is one-dimensional, then the
quasiunipotence condition means that the monodromy operators along
any sufficiently small loop around each singular point are
quasiunipotent. This condition can be effectively verified if the
system is Fuchsian (has a first order pole) by inspection of the
spectrum of the corresponding residue matrix: all eigenvalues of
this matrix should be \emph{rational}.

However, a loop that encircles several singularities is not small,
thus quasiunipotence does not impose any conditions on the
corresponding ``large'' monodromy.
\end{Rem}

The system \eqref{pfs} (and the corresponding matrix 1-form
$\Omega$) is called \emph{regular} on $\mathbb M$, if for any
\emph{real analytic} path $\gamma\:(\R,0)\to(\mathbb M,\S)$ the
solution grows no faster than polynomially near the singular locus,
\begin{equation}\label{reg-cond}
    |X(\gamma(s))|^{\pm 1}\le c|s|^{-p}\qquad \forall s\in(\R^1,0),
\end{equation}
for some real numbers $c,p<+\infty$, eventually depending on the
path $\gamma$. Analyticity of the path $\gamma$ intends to rule out
spiraling along the singular locus. One can show that regularity is
sufficient to verify only along (real) line segments.

The following result can be considered as a ``\emph{removable
singularity assertion}'' for regular quasiunipotent systems.

\begin{Thm}[{\normalfont Kashiwara theorem
 \cite{kashiwara:qu}}]\label{thm:kashiwara}
A regular integrable system that is quasiunipotent at each point
outside an algebraic subset of codimension $\ge 2$, is globally
quasiunipotent.
\end{Thm}

\begin{Rem}[important]\label{rem:relax-quasiunipotent}
If in the definition of the quasiunipotence (local and global) we
replace the assumption on the spectrum to consist solely of roots of
unity by a \emph{weaker} assumption that all these eigenvalues have
modulus 1, then Theorem~\ref{thm:kashiwara} remains valid as well as
all other assertions, including our main result
(Theorem~\ref{thm:count} below). It is in this stronger form that
the results on general Fuchsian systems are formulated in
\cite{polyfuchs}. However, we are not aware of any naturally arising system
that satisfies this weaker assumption and is not quasiunipotent in the usual
(stronger) sense.
\end{Rem}

%
%
%
%
%

%
%
%

\begin{Thm}[{\normalfont folklore}]\label{thm:classic-pf}
The period matrix $X$ of Abelian integrals \eqref{X} satisfies an
\emph{integrable}, \emph{regular} and \emph{quasiunipotent} system
of equations of the form \eqref{pfs} on $\P^m$.
\end{Thm}

\begin{proof}
This is a classical ``well-known fact'' whose proofs are scattered
over a number of classical sources. We outline only the principal
arguments.

The period matrix $X$ described in \eqref{X}, is \emph{monodromic}:
the result of its continuation $\Delta_\gamma X$ along any closed
path $\gamma$ avoiding the ramification locus $\S_*$, differs from
the initial value $X$ by a constant monodromy matrix which describes
the parallel transport of the cycles $\delta_j$ along the path
\cite[Corollary 26.28]{iy:lade}.

This implies that the \emph{logarithmic derivative} $\Omega=\d
X\cdot X^{-1}$ is a single-valued matrix function defined on the
complement $\P^m\ssm\S$.

The growth of $X$ near the singular locus is at most polynomial:
indeed, the length of any cycle of the integration $\delta_j(\l)$ is
growing no faster then polynomially along any curve and so does the
integrand. Hence all entries of $\Omega$ grow no faster than
polynomially near $\S$. Being single-valued, $\Omega$ is necessarily
a \emph{rational} matrix 1-form. The integrability condition follows
immediately from the local representation of $\Omega$ as a logarithmic
derivative.

The properties of the monodromy group of $\Omega$ were studied in
great detail. The fact that the system is quasiunipotent was proved
by Brieskorn \cite{brieskorn} and Clemens \cite{clemens} by
completely different methods; the proofs were re-exposed a number of
times \cite{odo-2,Zol}.

The shortest way to prove the quasiunipotence is using the Kashiwara
theorem. A generic point of $\S$ corresponds to a simple normal
crossing of the curve $\G_\l$ or to a quadratic tangency of this
curve with the infinite line (which is a polar locus for $\omega$). In the first case the monodromy along
a small loop $\gamma$ around $a$ is described by the
\emph{Picard--Lefschetz formulas} \cite[\parasymbol 26I]{iy:lade}:
the corresponding monodromy matrix $M_\gamma$ has a Jordan block of
size $2\times 2$ with the eigenvalue $1$. In the second case $M_\gamma^2=1$, hence all eigenvalues of $M_\gamma$ are necessarily $\pm1$. Both cases are clearly quasiunipotent. By the Kashiwara theorem, we conclude that $\Omega$
is globally quasiunipotent.
\end{proof}

\begin{Def}
Throughout this paper the Pfaffian system \eqref{pfs} satisfied by
the period matrix \eqref{X} will be referred to as \emph{the
Picard--Fuchs system}.
\end{Def}

\subsection{Polynomial norms}
When dealing with polynomials (both univariate and multivariate), we
will always use the $\ell^1$-norm.


\begin{Def}\label{def:norm}
The \emph{norm} of a multivariate polynomial
$P\in\C[z_1,\dots,\allowbreak z_n]$, $P(z)=\sum_\alpha c_\alpha
z^\alpha$, $c_\alpha\in\C$ (in the standard multiindex notation) is
the sum of absolute values of its coefficients, $\|P\|=\sum
|c_\alpha|$. Clearly, this norm is multiplicative,
\begin{equation}\label{m-norm}
    \|PQ\|\le\|P\|\cdot\| Q\|.
\end{equation}
\end{Def}

%

\subsection{Algebraic objects defined over $\Q$ and their
 complexity}\label{sec:overQ}
Complexification, replacing integrals over real ovals by the complex
analytic period matrix \eqref{X}, was one of the first tools of
investigation of Abelian integrals \cite{ilyashenko:elliptic}. In
this section we describe an opposite (in a sense) step and introduce
the class of differential equations defined over $\Q$, the subfield
of rational numbers in the field $\C$.

Speaking informally, an object (a polynomial, rational fraction,
variety, Pfaffian form, differential operator, semialgebraic set
\etc.) is defined over $\Q$, if it can be effectively constructed
from the ring $\Z[z_1,\dots,z_n]$ of \emph{lattice polynomials}. For
objects defined over $\Q$, besides the standard algebraic
characteristics like degree and dimension, one can always measure
its \emph{size}, roughly understood as the \emph{magnitude of the
integer numbers occurring in the explicit formulas} describing the
objects.

Perhaps the term ``size'' is not very successful, since the size of
a small reciprocal $1/n\in\Q$ with $n\in\N$ would be essentially equal to the
large number $n$. The term ``complexity'' would better suit our
purposes, but this term is too overloaded. As a compromise, the
reader may think of the size as (the exponent of) the \emph{bit}size
of the explicit representation of the given objects. The formal definitions follow.

\begin{Def}\label{def:size}
The \emph{size} $\size{R}$ of a lattice (integer) polynomial
$P\in\Z[z_1,\dots,z_n]$ is set to be equal to its norm, $\size
P=\|P\|$.

The \emph{size} of a rational fraction $R\in\Q(z_1,\dots,z_n)$ is
\begin{equation}\label{norm-R}
    \size{R}=\min_{P,Q}\{\|P\|+\|Q\|\:
    R={P}/{Q},\ P,Q\in\Z[z_1,\dots,z_n]\},
\end{equation}
the minimum being taken over all possible representations of $R=P/Q$
with $P,Q\in\Z[z_1,\dots,z_n]$.

%


The size of a (polynomial or rational) 1-form on $\P^m$ defined over
$\Q$, is the sum of sizes  of its coefficients in the standard
affine chart $\C^m\subset\P^m$.

The size of a vector or matrix rational function (resp., 1-form)
defined over $\Q$, is the sum of the sizes of its components.

A \emph{parametric family} of objects is defined over $\Q$, if it is
defined over $\Q$ on the product space $\P^{m-1}\times\P^1$
birationally equivalent to $\P^m$ (cf.~with
Remark~\ref{rem:pencil}). The size of the parametric family is
defined via the corresponding equivalence.

One can easily continue this series of definitions, extending it for
any class of algebraic objects and their parametric families
algebraically depending on auxiliary parameters. In the future we
will need to define the size of ordinary differential operators, see
\secref{sec:derivation}.
\end{Def}

\begin{Rem}
The size is associated not so much with an object, but rather with a
specific formula representing it. For instance, the polynomial
$1+t+t^2+\cdots+t^{n-1}$ of size $n$ in $\Z[t]$ can be represented
as a rational function $(t^n-1)/(t-1)$ of size only $4$ in $\Q(t)$.
Yet for our purposes this ambiguity will not be important, since we
will deal only with \emph{explicit upper bounds} for the size which
means construction of \emph{formulas} (representations) not
involving excessively large natural numbers.
\end{Rem}

The most important feature of the size of formulas defined over $\Q$
is its \emph{controlled increase in any algorithmically defined
construction}. For instance, the size of a sum and product of two
rational fractions of sizes $s_1,s_2$, can be estimated as follows,
\begin{equation}\label{compl-arithm}
    \size{\frac{p_1}{q_1}+\frac{p_2}{q_2}}
    \le\size{\frac{p_1q_2+p_2q_1}{q_1q_2}}\le
    3s_1s_2,\quad
    \size{\frac
    {p_1}{q_1}\cdot\frac{p_2}{q_2}}
    \le
    2s_1s_2.
\end{equation}
Composition of maps defined over $\Q$ is also an operation that
increases the size (complexity) in a limited way.

\begin{Ex}\label{ex:birational}
The projective space $\P^{m+n}$ is birationally equivalent to the
product $\P^{m}\times \P^n$, e.g., via the standard identification
$\C^{m+n}\simeq\C^m\times\C^n$ in the affine charts. Such
transformations do not result in a substantial change of the
complexity of objects, in particular, the above equivalence does not
change the complexity of rational functions defined on the
corresponding birationally equivalent varieties.
\end{Ex}

An example of effective complexity control is the following explicit
bound on solutions of linear systems of algebraic equations.

\begin{Ex}\label{ex:linsys}
Consider a parametric family of systems of linear algebraic
equations of the form
\begin{equation}\label{linsys}
    A(\l)z=b(\l),\qquad \l\in\P^m,\quad z=(z_1,\dots,z_n),
\end{equation}
and assume that this system is defined over $\Q$, that is, the
entries $a_{ij},b_j$ of the matrix $A$ and the right hand side
vector $b$ are elements from the field $\Q(\l)$ of known degrees not
exceeding $d$, and their size does not exceed $s$.

The system may well be non-solvable over $\Q(\l)$, but in case it is
known to have a solution, such a solution can always be found of
size not exceeding $O(n^n)\,s^n$.

Indeed, after eliminating all equations that are linearly dependent
over the field $\Q(\l)$, we can represent components of some
solution as ratios of suitable minors by the Cramer rule. Each of
these minors is the sum of at most $n!$ products of $n$ entries of
$A(\l)$, all of them in the field $\Q(\l)$, so that the degree is no
greater than $d^n$ and the size is no greater than $n!(6s)^n$ by
\eqref{compl-arithm}.
\end{Ex}

\subsection{``Quantization principle''}\label{sec:quantization}
The following  general principle is the primary reason why objects
defined over $\Q$ appear in a construction which initially has no
such structure: \emph{Any finite bound for objects defined over
$\Q$, is explicitly computable in terms of their size}. We give two
illustrations of this principle.

\begin{Ex}[algebraic, continuation of
Example~\ref{ex:linsys}]\label{ex:linsys2} A non-homogeneous system
$Az=b$ of linear algebraic equations defined over $\Q$ may have no
solution, but if the solution exists, there exists also a solution
of norm bounded in terms of the dimension $n$ and size $s$
(complexity) of the system, $|z|\le n!s^n$.

Indeed, without loss of generality we may assume that all entries of
the matrix $A$ and the free terms $b$ are \emph{integer} not
exceeding $s^n$ in the absolute value. Then, as explained in
Example~\ref{ex:linsys}, for some particular solution each component
can be represented as the ratio of some minors. The numerator does
not exceed $n!s^{n^2}$ (again using the Laplace expansion with
obvious estimates), while the denominator, being a nonzero
\emph{integer}, is no smaller than $1$.
\end{Ex}

\begin{Ex}[geometric]\label{ex:K-bound}
Assume that $K\subseteq\R^n$ is a basic semialgebraic set defined by
finitely many polynomial equalities and inequalities of the form
$p_\alpha(x)=0$ (resp., $p_\alpha\le 0$), where $p_\alpha$ are
polynomials defined over $\Q$ of degree $\le d$ and size no greater
than $s$.

The set $K$ may well be non-compact (e.g., a half-space), but if it
is known to be compact, its diameter can be explicitly bounded as
follows.

\begin{Thm}[{\normalfont\cite{roy.e.a.,basu-vorobjov}}]\label{thm:basu}
If a basic semialgebraic set
\begin{equation}\label{K-diam}
    K=\bigcap_\alpha \{x\:p_\alpha(x)\le0\},\quad p_\alpha\in\Z[x],\ \deg p_\alpha\le
    d,\ \|p_\alpha\|\le s,
\end{equation}
is bounded, then it belongs to the ball of radius $R$ centered at
the origin, with $R=s^{d^{O(n)}}$. The constant in $O(n)$ is
explicit.

The same result holds for semialgebraic sets defined by polynomials
from $\Q[x]$, if $s$ is the upper bound for their size
$\size{p_\alpha}$.
\end{Thm}

This claim, rather obvious for the one-dimensional case $n=1$, can
be proved for arbitrary $n$ by the algorithmic quantifier
elimination technique (corresponding to the projection of $K$ to
$\R^{n-1}\subset\R^n$). The quantifier elimination process can be
made constructive, which results in a controlled increase of the
complexity in each step.
\end{Ex}

Of course, the ``quantization principle'' is only a guiding line,
not a theorem; in each instance we will have to address a specific question
on effective bounds either directly (as in Example~\ref{ex:linsys})
or indirectly, using tools of effective real algebraic geometry (as
in Example~\ref{ex:K-bound}).

The ``quantization principle'' was already implicitly used in the
proof of the general result on meandering of trajectories of
arbitrary polynomial vector fields \cite{annalif-99}. In application
to linear systems this principle was introduced in
\cite{alexg:thesis,alexg:arxiv} and later in a more transparent and
general form in \cite{ongrig}.

\subsection{Counting zeros of multivalued matrix functions of several variables}
 \label{sec:count-complex}
The period matrix $X(\l)$, the  solution of the Pfaffian system
\eqref{pfs}, is ramified over the singular locus $\S$. We introduce
the \emph{counting function} which generalizes the number of
isolated zeros of functions of one variable to the class of
multivalued matrix-functions of several complex variables.

Let $\P^1\simeq L\subset\P^m$ be an arbitrary projective line
\emph{not entirely belonging} to the singular locus $\S$ of the
Pfaffian system \eqref{pfs} on $\mathbb M=\P^m$. The intersection
$L\cap\S$ in this case consists of finitely many isolated points,
and the restriction of the matrix function $X(\l)$ on $L$ will be
ramified over these points.

Let $T\subset L\ssm\S$ be an arbitrary \emph{triangle}, an open
domain bounded by three circular arcs eventually degenerating into
line segments or points. Since $T$ is simply connected, one can
unambiguously choose a continuous holomorphic branch of the matrix
function $X(t)=\bigl(x_{ij}(t)\bigr)_{i,j=1}^\ell$, $t\in T$.

Consider the \emph{linear span},
\begin{equation}\label{lin-span}
    \mathscr L_X(T)=\{f\in\mathscr O(T)\: f=\sum_{i,j=1}^\ell
    c_{ij}x_{ij}(t),\ c_{ij}\in\C\}.
\end{equation}
a finite-dimensional subspace in the space of functions of one
complex variable (recall that $T$ is a triangle in a complex
projective line), holomorphic in $T$.

Replacing the matrix function $X(\cdot)$ by a different solution
$X(\cdot)M$, $M\in\GL(\ell,\C)$ (in particular, by another branch of
analytic continuation of $X$), does not affect the subspace
$\mathscr L_X(T)$, thus the latter depends only on the Pfaffian
matrix 1-form $\Omega$.

We define the \emph{counting function} as the supremum
\begin{equation}\label{count-global}
    \cN(\Omega)=\sup_{T\subset\P^m\ssm\S}~\sup_{f\in \mathscr L_X(T)}\#\{t\in
    T\:f(t)=0\}\le +\infty,
\end{equation}
taken over all triangles $T$ disjoint with $\S$ and all nonzero
functions from $\mathscr L_X(T)$.

\begin{Rem}
The requirement that $T$ is a triangle is aimed at excluding simply
connecting domains spiraling around the singular locus. Easy
examples show that spiraling domains may contain arbitrarily large
number of isolated zeros of very simple functions.
\end{Rem}

\begin{Rem}\label{rem:count-in-domain}
Knowing the bounding function $\cN(\cdot)$, one can use
triangulation to estimate the number of isolated roots
$\cN(\Omega/U)$ of any linear combination in any polygonal domain
$U\subseteq L\ssm\S$ in any line $L\not\subset\S$. In particular,
the number of real zeros $\cN(\Omega/\Re L)$ can be at most
$d\cdot\cN(\Omega)$, where $d$ is the degree of $\Omega$. Indeed,
the intersection $L\cap\S$ consists of no more than $d$ points which
may subdivide the real (projective) line into to no more than $d$
intervals. Each interval lies inside a triangle $T$ free from points
of $\S$, hence the number of isolated zeros on it does not exceed
$\cN(\Omega)$.

Conversely, if there is a tiling $L\ssm\S=\bigcup \overline U_i$ by
simply connected polygonal domains $U_i$ and for each domain the
maximal number of zeros $\cN(\Omega/U_i)=\sup_{f\in\mathscr
L_X(U_i)}\#\{t\: f(t)=0\}$ is finite, then one can easily produce an
upper bound for $\cN(\Omega)$ by simply adding these bounds,
$\cN(\Omega)\le\sum_i\cN(\Omega/U_i)$.
\end{Rem}

\begin{Rem}[semicontinuity]\label{rem:genericity}
When counting zeros, one can easily pass from open to closed
polygons disjoint from $\S$, provided that the bound for the number
of zeros remains uniform.

Indeed, assuming that the number of zeros in \emph{any} closed
triangle $\ol T\subset L\ssm \S$ does not exceed some $N$, one can
immediately see that the same bound holds also for an arbitrary open
triangle $T$. If $T$ contains $N+1$ roots of some linear
combination, one can construct a closed triangle $\ol T\Subset T$
which contains all these roots, in contradiction with the initial
assumption.

By the same token, in the definition of the counting function
\eqref{count-global} it is sufficient to consider only closed
triangles $\ol T$ from a dense subset $\mathscr T'$ in the space of
all triangles $\mathscr T(\P^m\ssm\S)$ disjoint with $\S$. If
\begin{equation*}
    \sup_{\ol T\in\mathscr T'}~\sup_{f\in \mathscr L_X(T)}\#\{t\in
    T\:f(t)=0\}=N<+\infty,
\end{equation*}
then $\cN(\Omega)$ is also finite and equal to $N$. Indeed, if some
linear combination $f$ has $N+1$ roots in an ``excluded'' open
triangle $T\notin\mathscr T'$, then one can find an arbitrarily
close closed triangle $\ol T\in\mathscr T'$, eventually belonging to
a different line $L'\subset\P^m$, which contains at least $N+1$
roots counted with multiplicity in contradiction with the initial
assumption. This follows from the fact that isolated complex roots
of holomorphic functions cannot disappear by small perturbations by
the Weierstrass preparation theorem.
\end{Rem}

\subsection{Main theorem in the abstract form and discussion}
The constructive solution of the infinitesimal Hilbert 16th problem
(Theorem~\ref{thm:main2}) is obtained as a corollary to the
following general fact about solutions of systems of differential
equations.


\begin{Thm}[{\normalfont existential bound}]\label{thm:exist}
Let $\Omega$ be a rational matrix 1-form of degree $d$ on the
projective space $\mathbb M=\P^m$, and \eqref{pfs} the corresponding
linear system of size $\ell\times\ell$. Assume that\textup:
\begin{enumerate}
 \item[(I)~] $\Omega$ is integrable\textup;
 \item[(R)] $\Omega$ is regular\textup;
 \item[(U)] $\Omega$ is quasiunipotent.
\end{enumerate}
\noindent Then the value of the counting function $\cN(\Omega)$ is
finite,
\begin{equation*}
    \cN(\Omega)<+\infty.
\end{equation*}
\end{Thm}

\begin{Thm}[{\normalfont constructive bound}]\label{thm:count}
In the assumptions of Theorem~\ref{thm:exist} and the additional
assumption,
\begin{enumerate}
  \item[(Q)] $\Omega$ is defined over $\Q$ and its size is
  $s=\size{\Omega}$,
\end{enumerate}
the above finiteness is explicit\textup:
\begin{equation}\label{effective-bound}
    \cN(\Omega) \le s^{2^{\poly(d,\ell,m)}}.
\end{equation}
Here $\poly(d,\ell,m)\le O^+\bigl((d\ell^4m)^5\bigr)$.
\end{Thm}

Recall that the $O^+(\bullet)$-notation was introduced in
Remark~\ref{rem:oplus}. The reduction from
Theorem~\ref{thm:count} to Theorem~\ref{thm:main2} is made
possible by virtue of the following observation which improves
the ``folklore'' Theorem~\ref{thm:classic-pf}.

\begin{Thm}[{\normalfont constructive derivation}]\label{thm:pf-compl}
The logarithmic derivative $\Omega=\d X\cdot X^{-1}$ of the period
matrix for Abelian integrals \eqref{X}, and hence the corresponding
Picard--Fuchs system \eqref{pfs} is defined over $\Q$.

The size $s=\size{\Omega}$, dimension $\ell$ and the degree $d=\deg
\Omega$ of the corresponding rational matrix function are explicitly
bounded from above as follows,
\begin{equation}\label{bounds-pfs}
    s\le 2^{\poly(n)},\quad d\le O(n^2),\quad m\le O(n^2),\quad \ell= n^2,
\end{equation}
where $n+1$ is the degree of the Hamiltonians.
\end{Thm}

\medskip
\emph{Proof of Theorem~\ref{thm:main2}}.\quad Plugging the
estimates \eqref{bounds-pfs} into the bound
\eqref{effective-bound}, we obtain the bound for the number of
zeros of Abelian integrals,
\begin{equation}\label{oplus-main}
    \cN(n,n)\le 2^{\poly(n)\cdot\, 2^{O\(n^{(2+8+2)\cdot 5}\)}}
    \le 2^{2^{O(n^{60}\log n)}}\le 2^{2^{O^+(n^{60})}}.
\end{equation}
It remains to notice that $O^+(n^{60})\le O(n^{61})$. This
calculation illustrates the need for the $O^+$-notation, as
similar estimates will appear in the future.
\medskip



%
%
%
%

In fact, Theorem~\ref{thm:exist} can be relatively easily proved by
the application of tools from the Fewnomial theory developed by
A.~Khovanski\u\i~\cite{asik:fewnomials} and finiteness results for
analytic functions achieved by A.~Gabrielov
\cite{gabrielov:finiteness}. Unfortunately, this straightforward
approach does not allow for application of the ``quantization
principle'' mentioned in \secref{sec:quantization}.

We give an alternative proof of Theorem~\ref{thm:exist} that admits
the required ``quantization'' and ultimately yields
Theorem~\ref{thm:count}. The main ideas of this proof are outlined
below in \secref{sec:strategy}.

The conditions of Theorem~\ref{thm:count} are very close to optimal.
Indeed, without the integrability assumption the system has no
well-defined solutions. Omission of the regularity assumption allows
to construct a linear system on $\P^1$ (i.e., in the smallest
dimension), with infinitely many zeros of solutions on a real
interval accumulating to a singular point, see \cite{montreal}.

The assumption of quasiunipotence also cannot be considerably
relaxed beyond the limits indicated in
Remark~\ref{rem:relax-quasiunipotent}: without the assumption zeros
of solutions also can accumulate to a singular point along the real
line \cite{montreal}.


%

\subsection{Strategy of the proof}\label{sec:strategy}\begin{small}
In this section we briefly outline the strategy of the proof of
Theorems~\ref{thm:exist}, \ref{thm:count} and~\ref{thm:pf-compl}.

\subsubsection{The analytic core: de la Vall\'ee Poussin theorem and its generalizations.}
The basic tool for the estimates of the number of isolated zeros is
a complex analog of the classical de la Vall\'ee Poussin theorem
\cite{poussin}. This theorem asserts that for a homogeneous
\emph{monic} linear ordinary differential equation with holomorphic
coefficients
\begin{equation*}
    y^{(k)}+a_1(t)\,y^{(k-1)}+\cdots+a_k(t)\,y=0, \qquad
    t\in\gamma\subset\C,
\end{equation*}
the variation of argument of any solution $y(t)$ along a circular
arc $\gamma$ of known length can be explicitly bounded in terms of
the \emph{uniform upper bounds} $A_i=\sup_{t\in\gamma}|a_i(t)|$,
$i=1,\dots,k$, of the non-principal coefficients of this equation
along the arc \cite[Theorem 2.6, Corollary 2.7]{fields}. This
property of high order differential equations is in stark contrast
with \emph{systems} of first order linear equations with bounded
coefficients, as was discovered in \cite{mit:counterexample}. The
assumption that the equation is monic (i.e., its principal
coefficient is identically equal to $1$) is not an obstruction as
long as the arc does not pass near singular points of the equation,
which correspond to zeros of the leading coefficient of a general
homogeneous equation
\begin{equation}\label{pre-lode}
    a_0(t)\,y^{(k)}+\cdots+a_k(t)\,y=0.
\end{equation}
Computing variation of argument of solutions along an arc that
passes through (or very close to) a singular point (a root of
$a_0(\bullet)$) is impossible in general.

\subsubsection{From a Pfaffian integrable system to an isomonodromic
 family of homogeneous linear ordinary differential
 equations.}\label{sec:polyfuchs}
The system of Pfaffian equations \eqref{pfs} can be reduced to a
scalar equation of order $k\le \ell^2$ in two steps. First, we note
that the phase space $\mathbb M=\P^m$ of the system \eqref{pfs} is
birationally equivalent to the product $\P^{m-1}\times\P^1$,
cf.~with Remark~\ref{rem:pencil}; moreover, this equivalence can be
chosen so that any given projective line $\P^1\simeq L\subset\P^m$
becomes a member of the pencil of projective lines. Denoting the
parameters of the pencil by $\lbar\in\P^{m-1}$ and the corresponding
lines by $\P^1_\lbar=\{\lbar\}\times\P^1$, we can restrict the
Pfaffian systems on the lines from this pencil to obtain a
parametric family of Pfaffian equations on the projective line
$\P^1$ as in Remark~\ref{rem:pencil}. Since the latter space is
one-dimensional, choosing an affine chart $t$ on $\P^1$ allows to
re-write \eqref{pfs} as a linear system of first order ordinary
differential equations. Namely, the matrix Pfaffian 1-form
restricted on each line $\P^1_{\lbar}=\{\lbar\}\times\P^1$ in the
fixed chart $t$ takes the form
\begin{equation}\label{pfs-rest}
    \Omega|_{\{\lbar\}\times\P^1}=A(\lbar)\d t,\qquad
    \lbar\in\P^{m-1}.
\end{equation}
The system \eqref{pfs} becomes in this chart a system of linear
equations
\begin{equation}\label{pfs-odo}
    \frac{\d X}{\d t}=A(\lbar)\cdot X,\qquad X=X(\lbar,t),
\end{equation}
with the singular loci $\S_{\lbar}=\S\cap\P^1_{\lbar}$ (the matrix
$A(\lbar,\bullet)$ is not defined if the entire line $\P^1_{\lbar}$
belongs to $\S$, yet such values of the parameter constitute a
proper semialgebraic set in $\P^{m-1}$). Clearly, the regularity of the
initial system \eqref{pfs} implies the regularity of all systems in the
family \eqref{pfs-odo}.

The condition of integrability of the initial system (either on
$\P^m$ or on $\P^{m-1}\times\P^1$, as the two are equivalent)
implies that the family of the equations \eqref{pfs-odo} is
\emph{isomonodromic} in the following sense.

Let $\gamma\subset\P^1_{\lbar_0}$ be an arbitrary closed path in the
projective line, disjoint with the singular locus $\S_{\lbar}$. Then
by continuity there exists a small neighborhood $U$ of $\lbar_0$ in
$\P^{m-1}$ such that for all values of the parameter from this
neighborhood, the corresponding singular loci $\S_{\lbar}$ are still
disjoint with $\gamma$, and hence for all such $\lbar$ the monodromy
of the system \eqref{pfs-odo} along $\gamma$ is still well defined.
The isomonodromy condition means that the corresponding operators
$M_{\gamma,\lbar}$ do not depend on $\lbar\in U$, or, more
precisely, remain in the same conjugacy class.

This condition will be especially important when the singular locus
$\S_{\lbar}$ undergoes a ``bifurcation'' at $\lbar_0$, e.g.,
$\S_{\lbar_0}$ contains an isolated singularity of high multiplicity
at $t_0\in\C$ (a pole of order $k\ge 2$ for the corresponding matrix
$A(\lbar_0,\cdot)$), while all nearby matrices have simple
singularities (poles) at $k$ nearby points $t_i(\lbar)$,
$i=1,\dots,k$. The monodromy around the circle encompassing all
these points (the so called \emph{classical monodromy} of a
singularity, \cite{odo-2}) is conjugated to the monodromy of a small
loop around $t_0$.

\begin{Rem}
From now on we work only with the product space $\P^{m-1}\times\P^1$
and parametric equations and systems of equations on this product.
To simplify the notation, we replace the parameter space $\P^{m-1}$
by $\P^m$ and denote the coordinates on it by $\l$ instead of
$\lbar$.
\end{Rem}

On the second step we reduce the parametric system of linear
ordinary differential equations \eqref{pfs-odo} to a family of
scalar high order equations in the most straightforward way, by
successive derivations and linear elimination. As a result, we
obtain a linear equation with coefficients that are polynomials in
$t$ and the parameters $\l\in\P^m$,
\begin{equation}\label{dereq}
    a_0(\l,t)\,y^{(k)}+\cdots+a_k(\l,t)\,y=0,\qquad
    t\in\C,\ \l\in \C^m.
\end{equation}
Integrability and regularity of the initial system means that each
equation in the family \eqref{dereq} is Fuchsian and the family as a
whole is isomonodromic and quasiunipotent. In what follows we call
this family the \emph{derived equation}(s). However, two problems
arise in connection with this process:
\begin{enumerate}
 \item[(1)] The leading coefficient of the derived equation \eqref{dereq}
has isolated zeros, in general much more numerous than the singular
points of the original system.
 \item[(2)] For specific values of the parameters belonging to a proper
algebraic subset $\scS\subset\P^m$, the leading coefficient may
degenerate identically, $a_0(\l,\bullet)\equiv0$ for
$\l\in\scS$.
\end{enumerate}
The second problem implies that as $\l\to\scS$, the equation
undergoes what is usually referred to as a \emph{singular
perturbation}, the situation when the coefficient of the highest
order derivative tends to zero. Behavior of solutions of singularly
perturbed equations may be extremely complicated, and this scenario,
if it indeed occurs, renders the entire approach via de la Vall\'ee
Poussin's theorem unworkable, since after reducing to the monic form
the coefficients of the corresponding linear equations would be
large on the entire plane (or most of it), not just near singular
points. Note that we can ignore the exceptional value of parameters
$\l\in\scS$ itself by virtue of the Remark~\ref{rem:genericity},
provided that the bound for the number of zeros remains uniform.

Somewhat miraculously, the ``singular perturbation'' that occurs in
the reduction of a \emph{regular} system to the derived equation, is
only apparent: together with the leading coefficient, all other
coefficients of the equation \eqref{dereq} necessarily vanish for
the same values of the parameter $\scS$. This phenomenon was first
discovered (in a simpler context) by A.~Grigoriev
\cite{alexg:thesis,alexg:arxiv}.

The analysis carried out in \cite{polyfuchs} is  reproduced and
generalized in \secref{sec:derivation} to show that if we consider
the norms of these coefficients $\|a_i(\l,\bullet)\|$ (in the sense
of Definition~\ref{def:norm}), then the maximum of the ratios of
these norms, called the \emph{slope} of the linear homogeneous
equation, is uniformly bounded,
\begin{equation}\label{pre-slope}
    \max_{i=1,\dots,k}\sup_{\l\in\P^m\ssm\scS}
    \frac{\|a_i(\l,\bullet)\|}{\|a_0(\l,\bullet)\|}<+\infty.
\end{equation}

We note that the finiteness of this uniform  bound depends crucially
on the \emph{regularity} of the original system, and does not hold
in the more general context of \cite{polyfuchs}.

Finiteness of this supremum for a general regular family
\eqref{pfs-odo} implies, by the ``quantization principle'', an
explicit computable bound provided that the initial system
\eqref{pfs-odo} is defined over $\Q$.

\subsubsection{Invariant slope.}\label{sec:pre-slope}
The constructions described in \secref{sec:polyfuchs}, along with the techniques of \cite{polyfuchs}, provide a uniform bound for the number of zeros of solutions as long as the singular points remain well apart. A substantial difficulty
which needs to be addressed (and cannot be circumvented for general
Fuchsian systems) is the study of colliding singular points.

However, under the regularity and integrability assumptions, this problem can be resolved. It turns out that for \emph{Fuchsian equations} (operators) the slope as it is defined in \eqref{pre-slope} remains uniformly bounded from above even when one is allowed to replace the original affine chart by an arbitrary different affine (and more generally, even a conformally equivalent) chart on $\P^1$. This fact,
somewhat surprising (it seems to have gone unnoticed until now),
means that the collision of singularities can be treated by a proper
resolution of singularities, as explained in \secref{sec:pre-slits}.

In Definition~\ref{def:inv-slope} we introduce the notion of the
\emph{invariant slope} of a linear operator, and subsequently show
that it is finite and uniformly bounded. For technical reasons we
need to consider not only conformal changes of the independent
variable, but also symmetrization of differential operators with
respect to an arbitrary circle or line in $\P^1$. The corresponding
inequality is asserted in Principal Lemma~\ref{plem:slope},  whose
proof is presented in section \secref{sec:derivation}.

%

\subsubsection{Scale invariance and construction of an admissible configuration of
slits.}\label{sec:pre-slits} The bounds on the slope discussed in
\secref{sec:pre-slope},  imply that without loss of generality one
may assume that the leading coefficient of the derived equation
\eqref{dereq} is of unit norm with the remaining (non-leading)
coefficients uniformly bounded. By the de la Vall\'ee Poussin
theorem, this means that the variation of argument of any solution
can be explicitly majorized along any arc (circular or rectilinear),
which does not pass near the singular points where the leading
coefficient vanishes. The meaning of this dangerous proximity can be
made precise using the scale invariance of the invariant slope: the
upper bounds for the variation of arguments of an arbitrary solution
of the differential equation are possible for \emph{arcs whose
length is not very large relative to the distance to the singular
locus}. This scale invariance is one of the the key tools in the
subsequent construction. For brevity we refer to such arcs as
``short arcs''. For instance, any sufficiently small circular arc
centered at an isolated singular point, is ``short'' in this sense.

One may attempt to slit the plane with deleted singular points of
the derived equation by such ``short arcs'' into finitely many
simply connected domains $U_i$ and apply to each domain the argument
principle. This would  imply an explicit upper bound for the
counting function $\cN(\Omega/U_i)$, see
Remark~\ref{rem:count-in-domain}. Unfortunately, such simplistic
solution is impossible, since any arc with an endpoint at a singular
point is necessarily ``long''.

To resolve this problem, we show that one can bound  the number of
zeros of (multivalued) solutions of the derived equation in
punctured disks around singular points, and more generally in
annuli, under the assumption that the monodromy of this equation
along the equator of the annuli has eigenvalues of unit modulus
only, and that the bounding circles are ``short''. Our approach goes
back to the work \cite{mrl-96} and is based on the idea called
\emph{the Petrov trick} after the pioneering work by G.~Petrov
\cite{petrov:trick}.

On the second step of the construction (in \secref{sec:renorm}) we
construct a system of ``short arcs'' subdividing the nonsingular set
into simply connected domains and \emph{annuli} bounded by nested
circles. This comes in the form of a suitable clustering: we need to
identify groups of singular points, such that distance between them
is much smaller compared to distances to other singular points
(outside the cluster). Then after a suitable ``magnification'' one
can treat points from the same cluster as ``being well apart''.
However, the construction needs to be iterated, since much smaller
clusters can be parts of the larger clusters. The accurate
construction involves ideas in the spirit of the Fulton-MacPherson
compactification of the configuration space (see \cite{fulton-macpherson}).

The main difficulty to deal with in this second step is an effective
construction of the system of slits so that all annuli that appear
at the end will have the required monodromy, and the explicit
calculation of the  ``normalized length'' of these slits. The source
of difficulty is, among other things, the apparent non-algebraicity
of the monodromy as a function of the parameters: in general, the
monodromy of solutions of linear systems cannot be written in closed
form. The way around this obstacle is to use the quasiunipotence and
integrability. The quasiunipotence guarantees that the monodromy is
quasiunipotent along the \emph{small loops} which may encircle
\emph{several} colliding singularities. The integrability
(isomonodromy) allows to conclude that the monodromy remains
quasiunipotent as long as the topological configuration of slits and
singular points remains unchanged. The latter condition is
topological (selection of a connected component in a suitable
configuration space), and it is known that all connected components
of semialgebraic sets are themselves semialgebraic and effectively
constructible \cite{roy.e.a.}. This allows the application of the
``quantization principle'', transforming the existential finite
bound for the ``normalized length'' of the admissible system of
slits into an explicit upper bound for systems originally defined
over $\Q$ in exactly the same way as was done in the first step.

Knowing the explicit length of admissible system of slits along
``short'' arcs (the ``cluster diameter'' as it is introduced in
Definition~\ref{def:cluster-diameter}) and the invariant slope of
the equations allows to complete the proof of
Theorem~\ref{thm:count}. This program is realized in
\secref{sec:renorm}.


\subsubsection{Effective derivation of the Picard--Fuchs system.}
To derive Theorem~\ref{thm:main2} from Theorem~\ref{thm:count}, we
need to show that the Picard--Fuchs system provided by
Theorem~\ref{thm:classic-pf}, is in fact defined over $\Q$ for a
suitably chosen parameter space. The arguments used in the ``proof''
of Theorem~\ref{thm:classic-pf}, as well as some other standard
approaches \cite{gavrilov:petrov-modules,gavrilov:annalif}, do not
allow to estimate the size (complexity) of the Picard--Fuchs system
\eqref{pfs}.

The necessary bounds follow from the explicit derivation of the
Picard--Fuchs system, suggested in \cite{bullscimat-2002} and based
on an earlier work \cite{redundant}. We reproduce it below in
Appendix~\ref{sec:compl-pf} and derive all the required complexity
estimates.\par
\end{small}

\section{From an integrable Pfaffian system to an isomonodromic
 family of Fuchsian linear equations}\label{sec:derivation}

In this section we work with an integrable rational Pfaffian system
\eqref{pfs} of dimension $\ell\times\ell$ and degree $d$ on the product
space $\P^m\times\P^1$, defined over $\Q$, of known complexity
(size) $s=\size\Omega$.

Because of the integrability, we may consider the system as an
isomonodromic family of linear systems on $\P^1$, parameterized by
the parameters $\l\in\P^m$, and write
\begin{equation}\label{pfs-fib}
    \d X|_{\P^1_\l}=\Omega_\l X,\qquad \Omega_\l=A(\l,t)\d t,\quad
    \l\in\P^m,
\end{equation}
in an affine chart $t$ on $\P^1_\l$. Denote the singular locus of
the system \eqref{pfs-fib} by $\S_\l\subset\P^1_\l\simeq\P^1$.

The main result proved in this section is the effective
transformation of the family of linear systems \eqref{pfs-fib} to a
parametric family of scalar equations of a high order,
\begin{equation}\label{famlodo}
    D_\l=\partial^k+R_1(\l,t)\,\partial^{k-1}+\cdots+R_k(\l,t),
    \qquad R_j\in\Q(\l,t),
\end{equation}
also defined over $\Q$, with an explicit control on the \emph{size} of
the family (defined as the sum of the sizes of all rational
coefficients $R_j\in\Q(\l,t)$, cf.~with Definition~\ref{def:size}).
For a family of regular operators this turns out to be sufficient for explicitly
controlling the \emph{slope} (see \eqref{pre-slope}) of the
operators \eqref{famlodo} \emph{uniformly over all admissible
$\l\in\P^m$}.

The rationale behind this step is the stark difference between
systems of first order equations and scalar high order equations in
what concerns zeros of their solutions, see \cite{mit:counterexample}. To
formulate the result more precisely, we need to elaborate the
definition of the slope from \eqref{pre-slope} and make it
conformally invariant.

%
%
%
%

\subsection{Differential operators of higher order and their affine slope}
We will work with \emph{homogeneous} linear ordinary differential
equations with rational coefficients in the \emph{monic} form
\begin{multline}\label{lode}
    y^{(k)}+R_1(t)\,y^{(k-1)}+\cdots+R_{k-1}(t)\,y'+R_k(t)y=0,\\
    R_1,\dots,R_k\in\C(t),
\end{multline}
(so that the leading coefficient is identically 1) and their
parametric families. Because of the homogeneity, the monic equation
\eqref{lode} can be re-written in the form $Dy=0$, where $D$ is a
differential operator
\begin{equation}\label{lodo}
    D=a_0(t)\,\partial^k+a_1(t)\,\partial^{k-1}+\cdots+a_{k-1}(t)\,\partial
    +a_0(t),\quad \partial=\frac{\d}{\d t},
\end{equation}
with \emph{polynomial coefficients} $a_0,a_1,\dots,a_k\in\C[t]$,
$a_0\not\equiv0$ (we denote symbolically this fact by writing
$D\in\C[\partial,t]$). Under the assumption
\begin{equation}\label{gcd}
    \gcd(a_0,\cdots,a_k)_{\C[t]}=1
\end{equation}
the coefficients $a_i$ are determined uniquely modulo a
\emph{scalar} common factor. This makes the following definition
self-consistent.

%

\begin{Def}\label{def:slope}
The (\emph{affine}) \emph{slope} $\angle D$ of a linear ordinary
differential operator $D\in\C[t,\partial]$ with polynomial
coefficients as in \eqref{lodo}, is the finite number
\begin{equation}\label{slope}
    \angle D=\max_{j=1,\dots,k}\frac{\|a_j\|}{\|a_0\|}<+\infty.
\end{equation}

The \emph{slope} of a homogeneous linear ordinary differential
equation with \emph{rational} coefficients \eqref{lode} is by
definition the slope of the linear operator $D$ with polynomial
coefficients \eqref{lodo}--\eqref{gcd}, such that \eqref{lode} is
equivalent to the equation $Dy=0$.
\end{Def}

The affine slope of an operator is a numeric measure of proximity of
the corresponding equation to the ``singular limit'', the result of
perturbing a linear equation so that the highest derivative enters
with a small parameter. Knowing the slope of a homogeneous equation
suffices to place an explicit upper bound for the variation of
argument of any nontrivial solution of this equation along an arc
that does not pass through the singular points of this equation.

\begin{Lem}[{\normalfont Lemma~8 from \cite{polyfuchs}}]\label{lem:var-arg-aff}
 Let $D$ be a differential operator \eqref{lodo} of order
$k$ with polynomial coefficients of degree $\le d$ and the slope
$S=\angle D$, and $\gamma$ a closed circular arc or line segment
disjoint with the singular locus $\S=\{a_0=0\}\subset\C$, which
belongs to the disk of radius $R$ centered at the origin.

Then the variation of argument of any nonzero solution of the
homogeneous equation $Dy=0$ along the arc $\gamma$ is explicitly
bounded,
\begin{equation}\label{var-arg-bd}
    \Var\Arg y(t)|_\gamma \le kS\,|\gamma|\,(R/r)^{O(d)}.
\end{equation}
where $|\gamma|$ is the length of the arc, $r=\dist(\gamma,\S)$.
\end{Lem}

\begin{Rem}\label{rem:inv-subspaces}
Homogeneous linear differential equations with rational coefficients
are the natural means of describing finite dimensional linear
subspaces of holomorphic functions on $\P$, invariant by monodromy
around a finite locus $\S$. For instance, if \eqref{pfs} is a
regular Pfaffian system on the projective line $\P^1$ with a
singular locus $\S$, then for any open set $U\subseteq\P\ssm\S$ the
linear space $\mathscr L_X(U)$ spanned by components of any
fundamental matrix solution $X$ of the system \eqref{pfs} is
invariant by the monodromy, and can be defined by the homogeneous
linear equation
\begin{equation}\label{wronskian}\def\pdt#1{\frac{\d^{#1}}{\d t^{#1}}}
    \det\begin{pmatrix}
    y&x_1(t)&\cdots&x_\ell(t)\\
    \pdt{} y&\pdt{} x_1(t)&\cdots&\pdt{} x_\ell(t)\\
    \vdots&\vdots&\ddots&\vdots\\
    \pdt\ell y&\pdt\ell x_1(t)&\cdots&\pdt\ell x_\ell(t)
    \end{pmatrix}
    =0,
\end{equation}
where $x_1(\bullet),\dots,x_\ell(\bullet)$ is a basis of $\mathscr
L_X(U)$. When expanded in the elements of the first column and
reduced to the monic form, the identity \eqref{wronskian} yields a
Fuchsian\footnote{Recall that a linear higher order differential
operator is Fuchsian if it is regular, i.e., its solutions,
multivalued functions on $\C\ssm\S_F$, exhibit polynomial growth
as in \eqref{reg-cond}. As is well-known, Fuchsian operators admit
explicit description in terms of the order of zeros of their
coefficients at the singular locus \cite[Proposition
19.18]{iy:lade}.}  differential operator of order $\ell$ with
rational coefficients, provided that the functions $x_i(t)$ have
moderate growth near all points of the singular locus $\S$
\cite[Proposition 19.19]{iy:lade}.

This observation allows to define the slope of
\emph{any} finite-dimensional subspace $\mathscr L\subset\mathscr
O(U)$, $U\subseteq\P\ssm\S$, invariant by the monodromy
transformations associated with the fundamental group
$\pi_1(\P\ssm\S, t_0)$, $t_0\in U$, assuming that functions from $\mathscr L$ grow moderately near $\Sigma$. The slope $\angle\mathscr L$ is then defined as the slope of the
corresponding differential operator $D=D_{\mathscr L}$, the
differential operator with rational coefficients of the minimal
order $\ord D=\dim_\C\mathscr L$, which vanishes on $\mathscr L$:
\begin{equation}\label{slope-subspace}
\begin{gathered}
    \angle \mathscr L=\angle D_{\mathscr L},\qquad
    D_{\mathscr L}\in\C[t,\partial],\\
    \ord D_{\mathscr L}=\dim_\C\mathscr L,
    \qquad \forall f\in\mathscr L\quad D_{\mathscr L}f=0.
    \end{gathered}
\end{equation}
Note that this \emph{does not depend} on the choice of operator $D_{\mathscr
L}$, as any two linear differential operators with identical kernels agree
up to multiplication of the coefficients by a common factor. The slope \emph{does}
however depend on the choice of the affine chart $t$.
\end{Rem}

\begin{Rem}\label{rem:subspace-unique}
The exact choice of the simply connected domain $U$ is not important
as long as it is open and disjoint with the singular locus $\S$,
since the rational coefficients of the differential equation are
uniquely determined by their values in any open subset of $\P$. This
allows us to omit the indication of the domain in the notations.
\end{Rem}

\subsection{Conformal invariance and symmetrization}
The notion of a slope as it is defined in \eqref{slope} and
\eqref{slope-subspace} suffers from several drawbacks, the
most serious being its dependence on the chart with respect to
which the norms of the polynomial coefficients are computed.
Applications of different conformal automorphisms (changes of
the independent variable) of the form
\begin{equation}\label{conf-aut}
    \f\:t\longmapsto\frac{\alpha t+\beta}{\gamma t+\delta},
    \qquad\det\begin{pmatrix}\alpha&\beta\\\gamma&\delta\end{pmatrix}\ne0,
\end{equation}
transform any linear subspace $\mathscr L\subset\mathscr O(U)$ into
another subspace
\begin{equation}\label{action-aut}
    \f^*\mathscr L=\{\f^* f\:f\in\mathscr L\}\subset\mathscr
    O(\f^{-1}(U)),\qquad \f^*f=f\circ\f.
\end{equation}
having the same dimension and invariant by the monodromy operators
around the transformed singular locus $\S'=\f^{-1}(\S)$. We will
modify the definition of the slope so that it will become invariant
under the actions  of the conformal isomorphisms \eqref{action-aut}.

Besides the action of conformal isomorphisms, we will need yet
another operator on linear spaces of functions, the
\emph{symmetrization} (or \emph{conjugacy}) with respect to a circle/line. The
need for symmetrization will become clear in the context of Lemma~\ref{lem:annulus-count}. Recall that by the Schwarz
symmetry principle, for any function $f\in\mathscr O(U)$ holomorphic
in a domain $U\subseteq\C$, the function $f^\dag$ defined in the
domain $U^\dag$ by the formula
\begin{equation}\label{reflect}
    f^\dag(t)=\overline{f(\bar t)},\qquad U^\dag=\{\bar t\: t\in
    U\},
\end{equation}
is also holomorphic. We will refer to $f^\dag$ as the
\emph{reflection} of $f$ in the real axis $\R$.

If $U=U^\dag$ is itself a domain symmetric with respect to $\R$, then a
$\C$-linear subspace $\mathscr L\subset\mathscr O(U)$ will be called
\emph{real} (on $\R$), if its reflection $\mathscr
L^\dag=\{f^\dag\:f\in\mathscr L\}$ coincides with $\mathscr L$. A
finite-dimensional real subspace admits a basis (over $\C$) of
functions, real (i.e., taking real values) on $U\cap\R$: it is
sufficient to consider functions of the form $\frac12(f+f^\dag)$.

Any linear subspace $\mathscr L\subset\mathscr O(U)$ can be
\emph{symmetrized} as follows,
\begin{equation}\label{symmetrization}
    \mathscr L^\ominus=(\mathscr L+\mathscr L^\dag)|_{U\cap U^\dag},
\end{equation}
(note that the functions from $\mathscr L^\ominus$ need to be
restricted on the \emph{symmetrization} $U^\ominus=U\cap U^\dag$).

The role of the real axis can be played by any circle or real line
$\gamma\subset\P$ equivalent to the ``standard'' real
axis $\R\subset\P$ by a conformal automorphism $\f\in\Aut(\P^1)$. The conformal equivalence will be denoted by the
relation $\gamma\conf\R$. If $\f$ is an automorphism
which transforms $\gamma$ to the real axis $\R$, then the
\emph{reflection in $\gamma$} is defined by the formula
\begin{equation}\label{reflection-gamma}
    f^\dag(\f(t))=\overline{f(\overline{\f(t)})}.
\end{equation}
Modifying all constructions above accordingly, we arrive at the
notion of a \emph{symmetrization} of a linear subspace of
functions $\mathscr L$. Such symmetrization will be denoted by
$\mathscr L^\ominus_\gamma$. An subspace $\mathscr L$ such that
$\mathscr L^\ominus_\gamma=\mathscr L$ is called \emph{real on
$\gamma$}.

The dimension $\dim_\C\mathscr L^\ominus_\gamma$ of the
symmetrization depends, in general, on the arc $\gamma$. For
instance, if a space $\mathscr L$ is a real on $\R$, then $\mathscr
L=\mathscr L^\ominus_\R$, hence its symmetrization with respect to
$\R$ has the same dimension, yet for a generic line
$\gamma_\e=e^{\iu \e}\R$, $0<\e\ll1$, arbitrarily close to $\R$, the
dimension of the symmetrization $\mathscr L^\ominus_{\gamma_\e}$
will be twice the dimension of $\mathscr L$.

\begin{Def}
The \emph{invariant slope} of a finite-dimensional monodromic
subspace $\mathscr L\subset\mathscr O(U)$, $U\subset\P\ssm\S$ a
simply connected polygon, is the supremum of slopes of all
symmetrized conformal images of $L$:
\begin{equation}\label{inv-slope-subs}
    \slope \mathscr L=\sup_{\f,\gamma}\angle (\f^*\mathscr L)^\ominus_\gamma,
    \qquad \f\in\Aut(\P^1),\ \gamma\conf\R.
\end{equation}
Here the supremum is taken  over all conformal automorphisms $f$ and
all circular arcs/lines $\gamma\subset\P$, conformally equivalent to
$\R$.
\end{Def}

By the natural duality between finite-dimensional monodromic
subspaces of holomorphic functions and linear ordinary differential
operators with rational coefficients (modulo a common factor),
described in Remark~\ref{rem:inv-subspaces}, the conformal
isomorphisms and symmetrization can be defined also on differential
operators with rational coefficients on $\P^1$. More precisely, for
a given operator $D$ of order $\ell$ with a singular locus $\S$, we
denote by $\f^*D$ the operator of order $\ell$ (unique modulo
multiplication by a rational function) whose null space is
$\f^*\mathscr L$, where $\mathscr L=\{Dy=0\}\subset\mathscr O(U)$ is
the null space of $D$ in any simply connected domain
$U\subset\P\ssm\S$. In the same way we denote by $D^\ominus_\gamma$
the \emph{symmetrization} of $D$ with respect to an arc
$\gamma\conf\R$, the operator of order $\le 2\ell$ with the null
space $\mathscr L^\ominus_\gamma$. As was already noted, the
specific choice of the domain $U$ is not important by
Remark~\ref{rem:subspace-unique}. An operator equal to its
symmetrization, $D=D^\ominus_\gamma$, will be called \emph{real on
$\gamma$}, though this does not mean in general that the
coefficients of this operator (in the monic representation) indeed
take real values on the ``axis of symmetry'' $\gamma$.

\begin{Def}\label{def:inv-slope}
The \emph{invariant slope} $\slope D$ of a linear ordinary
differential operator $D$ with rational coefficients on $\P^1$ is
the supremum
\begin{equation}\label{inv-slope}
    \slope D=\sup_{\f,\gamma}\angle
    (\f^*D)^\ominus_\gamma,\qquad \f\in\Aut(\P^1),\ \gamma\conf\R.
\end{equation}
\end{Def}

\begin{Rem}\begin{small}
The group $\Aut(\P^1)\simeq\operatorname{PGL}(2,\C)$ of conformal
automorphisms of the projective line is noncompact, therefore the
slope $\angle \f^*D$ may be unbounded as $\f$ varies over this
group. Similarly, the procedure of symmetrization may affect the
slope in an uncontrollable way.

\begin{Ex}
Let $\f_\mu\:t\mapsto\mu t$ be the linear rescaling map and $D\in\C[\partial]$
a linear operator with \emph{constant} coefficients. Then the slope
$\angle \f_\mu^*D$ is unbounded as $\mu$ varies over all nonzero
complex (and even real) numbers.
\end{Ex}

However, for \emph{Fuchsian} operators the supremum in
\eqref{inv-slope} is always finite.

\begin{Prop}\label{prop:inv-slope}
The invariant slope of any Fuchsian operator is finite.
\end{Prop}

We will give an indirect proof of this statement later, in
Remark~\ref{rem:inv-slope-finiteness}.
\par\end{small}
\end{Rem}


The invariant slope of an operator is a semialgebraic function of its
coefficients, yet it is very difficult to control: its value
requires division by quantities which can be arbitrarily small. Our
first main result circumvents this difficulty and shows that the
straightforward reduction of a parametric linear system
\eqref{pfs-fib} to a parametric high order equation \eqref{famlodo}
results in an explicitly bounded slope.

\begin{PLem}\label{plem:slope}\indent\par
{\textbf A}. Let $\Omega$ be a rational $(\ell\times\ell)$-matrix
Pfaffian system of degree $d$ on $\P^{m}\times\P^1$ with the
following properties \textup(cf.~with the assumptions of
Theorems~\ref{thm:exist},~\ref{thm:count}\textup),
\begin{enumerate}
 \item[(I)~] $\Omega$ is integrable\textup;
 \item[(R)] $\Omega$ is regular\textup.
\end{enumerate}

Then for any simply connected polygon $U\subset \P\ssm\S_\l$ the
linear spaces $\mathscr L_\l\subset\mathscr O(U)$ generated by all
components $x_{pq}(\l,t)$ of any fundamental matrix solution $X$ of
the system \eqref{pfs} in this domain, are defined by a family of
linear ordinary differential equations of the form $D_\l y=0$, where
$D=\{D_\l\}$ are Fuchsian operators as in \eqref{famlodo}.

The family $D$ does not depend on $U$ and satisfies the following
constraints\textup:
\begin{itemize}
 \item[(i)] the order $k=\ord D$ is no greater than $\ell^2$,
 \item[(ii)] the degree $\max_j\deg R_j$ is bounded by an explicit polynomial in $\ell$
 and $d=\deg\Omega$,
 \item[(iii)]the invariant slope $\slope D_\l$ of the operators is
uniformly bounded over all values of the parameter $\l\in\P^m$
outside a proper algebraic subset $\scS\subsetneq\P^m$.
\end{itemize}

{\textbf B.} Under the additional assumption
\begin{enumerate}
 \item[(Q)] $\Omega$ is defined over $\Q$ and its size \textup(complexity\textup) is
 $s=\size{\Omega}$,
\end{enumerate}
the bound for the invariant slope is explicit and double
exponential,
 \begin{equation}\label{slope-explicit}
 \forall\l\in\P^m\ssm\scS\qquad \slope D_\l\le s^{(d\ell)^{O(m)}}.
 \end{equation}
\end{PLem}

Note that in the proof of this result we did not use the
quasiunipotence assumption.

\begin{Def}
The family of equations $\{D_\l y=0\}$ constructed in Principal
Lemma~\ref{plem:slope}, will be referred to as the \emph{derived
equation(s)}.
\end{Def}

\subsection{Formal derivation}
In this section we recall a (fairly standard) reduction of the
parametric system \eqref{pfs-fib} to a parametric family of linear
ordinary differential equations, paying attention to the complexity
of the algorithm.

%
%

\begin{Lem}\label{lem:D-construction}
Under the assumptions \textup{(I),~(R)} of the Principal
Lemma~\ref{plem:slope} all components $y=x_{pq}(\l,t)$ of any
fundamental matrix solution $X$ of the system \eqref{pfs} on
$\mathbb M=\P^m\times\P^1$ satisfy a parametric linear ordinary
differential equation of the form $D_\l y=0$ as in \eqref{famlodo}.

The family $D=\{D_\l\}_{\l\in\P^m}$ satisfies the following
constraints\textup:
\begin{enumerate}
 \item[(1)] the order $k=\ord D$ is no greater than $\ell^2$,
 \item[(2)] the degrees $\deg R_j$ do not exceed $O(\ell^4 d)$.
\end{enumerate}
\noindent Under the additional assumption \textup{(Q)} of the
Principal Lemma~\ref{plem:slope},
\begin{enumerate}
 \item[(3)] the family $D$ is defined over $\Q$ and its
 size is bounded by a simple exponential,
 \begin{equation*}
    \size{D}\le s^{O(\ell^4 d)},\qquad s=\size{\Omega}.
 \end{equation*}
\end{enumerate}
\end{Lem}

\begin{proof}
This claim coincides (modulo notation) with Lemma~5 from
\cite{polyfuchs}, where one can find the accurate (albeit
transparent) estimates. To make the exposition self-contained, we
recall the main ideas of the proof.

 \medskip\begin{small}
The system \eqref{pfs} is defined on the product space
$\P^m\times\P^1$. Choosing an affine chart $t$ on the second factor,
we can consider it as a parametric family of linear differential
equations of the following form (cf.~with \eqref{pfs-fib}),
\begin{equation}\label{ls}
    \frac{\partial X}{\partial t}=A(\l,t)\,X.
\end{equation}

By induction, one can instantly see that the higher order
derivatives $\frac{\partial^j}{\partial t^j} X$ of the matrix $X$
which satisfies the linear system \eqref{ls}, satisfy the identities
\begin{equation}\label{deriv-induct}
    \frac{\partial^j X}{\partial t^j}=A_j(\l,t)\cdot X,\qquad
    A_{j+1}=\frac{\partial A_j}{\partial t}
    +A_j\cdot A,\qquad j=0,1,2,\dots,
\end{equation}
where $A_0=E$ is the identity matrix. Since $A$
is defined over $\Q$, all matrix functions in this sequence are also
defined over $\Q$.

Derivation and matrix multiplication result in a controlled growth
of the degrees and sizes of the rational matrix functions
$A_1,A_2,\dots$: the degree grows no faster than linearly in $j$,
$\deg A_j\le jd$, while the sizes $\size{A_j}$ grow no faster than
exponentially.

The rational matrix $(\ell\times\ell)$-functions over $\Q$ form
a linear space of dimension $\ell^2$ over the field $\Q(\l,t)$
of rational functions on $\P^m\times\P$. Thus after $k\le
\ell^2$ steps the matrices $A_j$ will exhibit a linear
dependence over this field of the form
\begin{equation}\label{lin-dep}
    A_k+R_1A_{k-1}+\cdots+R_kA_0=0,\qquad R_j\in\Q(\l,t).
\end{equation}
The identity \eqref{lin-dep} can be considered as a system of
$\ell^2$ linear algebraic equations over the field $\Q(\l,t)$.
Solutions of this system (after elimination of all linear
dependencies between equations) can be effectively computed using
the Cramer rule as ratios of suitable determinants formed by entries
of the matrices $A_1,\dots,A_k$. This allows to place an upper bound
for the sizes $\size{R_j}$ in terms of $s$ and the parameters
$d,\ell$.

The differential operator $D=\partial^k+\sum_1^k R_{k-j}\partial^j$,
$\partial=\frac{\partial}{\partial t}$ by construction vanishes on
$X$: $DX\equiv0$. Read componentwise, this matrix identity proves
the Lemma.\par\end{small}
\end{proof}

%

\subsection{Existential bounds for the slope of the derived family}
The family of differential operators $D=\{D_\l\}$ which is
constructed in Lemma~\ref{lem:D-construction}, can be always reduced
to a form with polynomial coefficients, which will be referred to as
the \emph{standard form},
\begin{equation}\label{stf}
\begin{gathered}
    D=p_0(\l,t)\,\partial^k+\cdots+p_{k-1}(\l,t)\,\partial+p_k(\l,t),
    \\
    p_j\in\Z[\l,t],\qquad \gcd\nolimits_{\Z[\l,t]}(p_0,\dots,p_k)=1,
\end{gathered}
\end{equation}
The operator $D_\l$ has order $k$ for almost all values of
$\l\in\P$. Moreover, by the Bertini--Sard theorem, the subset of the
parameter values $\l$,
\begin{multline}\label{gcd-l}
    \scS=\{\l\:p_0(\l,\bullet)=0\}\\\cup
    \{\l\:\deg\gcd\nolimits_{\C[t]}\bigl(p_0(\l,\bullet),\dots,p_k(\l,\bullet)\bigr)>0\}.
\end{multline}
is a proper algebraic subvariety, $\scS\subsetneq\P^m$, which is
nowhere dense in $\P^m$.

For all $\l\notin\scS$, denote by $\sigma(\l)=\angle D_\l$ the
\emph{affine} slope of the operator $D_\l$. By construction, it is a
\emph{semialgebraic function} on the parameter space.

\begin{Lem}\label{lem:exist-bd}
If the initial integrable system \eqref{ls} is regular, the affine
slope of the derived family $D_\l$ is globally bounded,
$$\sup_{\P^m\ssm\scS}\angle D_\l<+\infty.$$
\end{Lem}

\begin{proof}
We will prove that the function $\sigma\:\P^m\ssm\scS\to\R_+$ is
\emph{locally bounded} near each point $a\in\P^m$, i.e., that there
exists such neighborhood $U$ of $a$, such that $\sup_{\l\in
U\ssm\scS}\sigma(\l)<+\infty$. By compactness of $\P^m$, this would
imply that $\sigma$ is globally bounded. Clearly, it is sufficient
to consider only the points $a\in\scS$, i.e., the values of the
parameters $\l$ for which the leading coefficient of the operator
$D_\l$ vanishes identically: at all other points $\sigma$ is
continuous.

Moreover, we can assume without loss of generality that the
parameter $\l$ varies along a real analytic curve,
\begin{equation*}
    \l=\gamma(s),\quad s\in(\R^1,0);
    \qquad \gamma(s)\notin\scS\iff s\ne0.
\end{equation*}
Indeed, if the function $\sigma$ were not locally bounded, then
the point $(a,0)\in U\times\R^1_+$ would be the limit point for
the open semialgebraic set $S=\{(\l,z):\allowbreak
\l\notin\scS, \sigma(\l)>1/z\}$. By the curve selection lemma
\cite{milnor:sing-points}, the point $(a,0)$ can be reached
from inside $S$ by a real analytic curve
$s\mapsto(\gamma(s),z(s))$, which means that the function
$\sigma(\l(s))$ grows to infinity as $s\to0$.

Thus we need to consider only the particular case of a differential
operator $D\in\mathscr O(\R^1,0)[t,\partial]$ depending on a single parameter,
\begin{equation}\label{D-s}
    D=p_0(s,t)\,\partial^k+\cdots+p_k(s,t),\qquad s\in(\R^1,0),\ t\in\C,
    \ \partial=\tfrac{\partial}{\partial t},
\end{equation}
with coefficients $p_j$ polynomial in $t$, real analytic in
$s\in(\R^1,0)$ and having no common factor for all $s\ne 0$. Because
of the real analyticity, we can complexify \eqref{D-s} to become a family
of differential operators $D_s=D|_{s=\const}$ defined for all sufficiently
small complex $s\in(\C^1,0)$.

The singular locus of this family is the analytic curve
$\{(s,t)\in(\C,0)\times\C\:p_0(s,t)=0\}$. Apart from the axis
$\{s=0\}$ corresponding to the identically vanishing leading
coefficient, each axis $\{s\}\times\C$ intersects this curve by
finitely many points corresponding to singularities of the equation
$D_{\l_s}y=0$. These points lie on branches of the above analytic
curve, therefore one can always find a disk $\mathbb D\subset\C$ of
radius $1$, such that the product $(\C,0)\times\mathbb D$ intersects
the singular locus only by the disk $\{0\}\times\mathbb D$.

Consider a fundamental system of solutions $x_1(s,t),\dots,x_k(s,t)$
of the equation $D_s y=0$ in the product $(\C^1,0)\times\mathbb D$,
which consists of the linearly independent components of a matrix
solution $X(\l(s),t)$ of the initial system \eqref{ls}. By the
choice of $\mathbb D$, these functions are holomorphic outside
$\{s=0\}$ and linearly independent, but may well be ramified with
a nontrivial monodromy $\Delta$ corresponding to a loop around the
axis $\{s=0\}$,
\begin{equation*}
    \Delta (x_1,\dots,x_k)=(x_1,\dots,x_k)\cdot M,\qquad M\in\GL(k,\C).
\end{equation*}
Fix a constant matrix $A$ such that $\exp2\pi\iu A=M$ (such a
matrix always exists, since $M$ is invertible). Then the tuple of
functions
\begin{equation*}
    (x_1',\dots,x_k')=(x_1,\dots,x_k)\cdot s^A
\end{equation*}
is single-valued in $(\C^1,0)\times\mathbb D$. Because of the
regularity of the initial system the new tuple of functions has
at most poles of finite order on the axis $\{s=0\}$. After
replacing $A$ by $A+qE$ for sufficiently large $q\in\N$, we
construct a tuple of functions $x_j'(s,t)$, still forming a
fundamental system of solutions for the family $\{D_s\}$ in
$(\C,0)\times\mathbb D$ for all $s\ne0$, such that $x'_j$ are
\emph{holomorphic} on the axis $\{s=0\}$. Note that the
restrictions of these functions on the axis itself may well
become degenerate (linearly dependent, e.g., identical zeros).

Consider the $\C$-linear subspaces $\mathscr
L_s\subseteq\mathscr (\mathbb D)$ spanned by the functions
$x_1'(\bullet,s),\dots,x_k'(\bullet,s)$. These subspaces depend
holomorphically (in the natural sense) on $s\in(\C^1,0)$ as
long as $s\ne 0$. The above mentioned degeneracy theoretically
means that the analyticity breaks down at $s=0$.

One of the keystone results, Lemma~7 from \cite{ongrig}
(cf.~with Proposition~18.18 from \cite{iy:lade}), contends that
this is not the case, and \emph{the application
$s\mapsto\mathscr L_s$ is a holomorphic curve in the
``Grassmanian'' of $k$-dimensional subspaces in the Banach
space $\mathscr O(\mathbb D)$.} This result can be seen as a
removable singularity-type theorem. One can avoid technical
difficulties of dealing with infinite-dimensional Grassmanians
by stating that \emph{one can choose different bases in the
subspaces $\mathscr L_s$ which would depend analytically on $s$
for all $s\in(\C^1,0)$ and remain linear independent as $s=0$}.

\begin{Lem}[{\normalfont Lemma~7 from \cite{ongrig}}]\label{lem:7-ongrig}
For any collection of holomorphic functions
$x_1',\dots,x_k'\in\mathscr O\bigl((\C^1,0)\times\mathbb
D\bigr)$ such that $x_j'(s,\bullet)$ are linear independent in
$\mathscr O(\mathbb D)$ for all $s\ne0$, one can construct a
meromorphic matrix function $R(s)$, nondegenerate and
holomorphic for $s\ne 0$, such that the tuple of functions
$f_j=f_j(s,t)$,
\begin{equation*}
    (f_1,\dots,f_k)=(x_1',\dots,x_k')\cdot R(s)
\end{equation*}
is holomorphic in $(s,t)$ and linearly independent on each
fiber $\{s=\const\}$, \emph{including the axis} $\{s=0\}$. \qed
\end{Lem}

By construction, $D_sf_j=0$ for all $s\in(\C^1,0)$ and all
$j=1,\dots,k$.

In a standard way, one can construct a family of \emph{monic}
differential operators $L=\{L_s\}$ of order $k$ with
coefficients holomorphic in $(\C^1,0)\times \mathbb D$,
\begin{equation*}
\begin{gathered}
    L=\partial^k+q_1(s,t)\,\partial^{k-1}+\cdots+q_{k-1}(s,t)\,\partial+q_k(s,t),
    \\
    q_j\in\mathscr O\bigl((\C^1,0)\times\mathbb D\bigr),\qquad j=1,\dots,k,
    \end{gathered}
    \end{equation*}
which is also annulled by the same tuple: $L_s f_j=0$ for all
$s\in(\C^1,0)$ and all $j=1,\dots,k$. Since the orders of the two
operators $L_s$ and $D_s$ are the same and the functions $f_j$ are
linearly independent on all fibers $\{s=\const\}$, we conclude that
the operators $D_s$ and $L_s$ must be proportional,
\begin{equation}\label{proportional}
    \frac{p_j(s,\bullet)}{p_0(s,\bullet)}=\frac{q_j(s,\bullet)}{1},
    \qquad\forall s\ne0,\quad \forall j=1,\dots,k.
\end{equation}
Since $q_j$ are holomorphic, this implies that the ratios in the
left hand side of \eqref{proportional} are holomorphic \emph{also}
on the axis $\{s=0\}$. We show that this implies the finiteness of the
ratios of the norms in $\C[t]$ as function of $s\to0$.

Let $\nu_j\ge 0$ be the vanishing orders of the polynomial
coefficients $p_j$ on the axis $\{s=0\}$: this means that
$p_j=s^{\nu_j}p_j'$, while $p_j'\in\mathscr O(\C^1,0)[t]$ and
$p_j'(0,\bullet)\ne0$. From \eqref{proportional} it follows
that $\nu_j\ge\nu_0$ for all $j=1,\dots,k$. Besides, since the
limit denominator $p_0'(0,\bullet)$ is nonzero, its norm
$\|p_0'(s,\bullet)\|_{\C[t]}$ is strictly positive. As a
result, we conclude that the slope $\sigma(s)=\angle D_s$ is
continuous at $s=0$:
\begin{equation*}
    \sigma(s)=\max_{j=1,\dots,k}
    s^{\nu_j-\nu_0}\cdot
    \frac{\|p_j'(s,\bullet)\|_{\C[t]}}{\|p_0'(s,\bullet)\|_{\C[t]}}.
\end{equation*}
This proves the local boundedness of the slope $\angle D_s$ along
the real analytic curve $\gamma$, concluding the proof of the Lemma.
\end{proof}

\subsection{Embedding in a conformally complete symmetric
 family}\label{sec:conf-compl}
The assertion of the Principal Lemma~\ref{plem:slope} concerns the
\emph{invariant} slope of the derived operator, whereas the
finiteness achieved in Lemma~\ref{lem:exist-bd} is established only
for the \emph{affine} slope. Besides, bounding the invariant slope
involves symmetrization, whose explicit construction may lead to an
uncontrollable growth of the slope, see \cite[Example 6]{polyfuchs}.

To fill the gap, we embed the initial parametric family of Pfaffian
systems \eqref{pfs-fib} into a larger family of systems of larger
dimension, whose derived equation will contain all symmetrizations of all
conformal transforms of the initial derived family. The uniform
bound for the affine slope in this new family gives a bound for the
invariant slope of the original one.

The embedding is rather straightforward in terms of the linear
spaces. Denote by $X(\l,t)$ the matrix solution of the system
\eqref{pfs} on $\P^m\times\P^1$ and
$\f\in\Aut(\P^1)\simeq\operatorname{PGL}(2,\C)\simeq\P^3\ssm\{\text{a
quadric}\}$ a variable conformal isomorphism. Then the components of
the matrix function $\boldsymbol X(\l,\f,\bullet)=X(\l,\f(\bullet))$
span for each admissible value value $(\l,\f)\in\P^m\times\P^3$ of
the parameters the linear space $\mathscr L_{\l,\f}=\f^*\mathscr
L_\l$ which is a conformal transform \eqref{action-aut} of the
linear space $\mathscr L_\l$ spanned by the components of the
initial matrix function $X(\l,\bullet)$. In other words, considered
as a parametric family with the parameter space $\P^m\times\P^3$,
the Pfaffian system with the matrix 1-form $\boldsymbol \Omega=\d
\boldsymbol X\cdot \boldsymbol X^{-1}$ contains all conformal
transforms of the initial system.

The dimension and the degree of the system $\boldsymbol \Omega$
obviously remain the same as that of the system $\Omega$. It is
almost as easy to see that the size (complexity) $\size{\boldsymbol
\Omega}$ differs from $\size{\Omega}$ by involving into a constant
power, $s\mapsto s^{O(1)}$. Indeed, the derivatives of $\boldsymbol
X$ with respect to the variables $\f\in\P^3$ can be expressed over
$\Q$ through entries of $\Omega$, $t$ and $\f$ using the chain rule,
and it remains to apply a few times the inequalities
\eqref{compl-arithm}. Clearly, the new family is also regular,
integrable and quasiunipotent.

It remains to embed the family of linear spaces $\mathscr L_{\l,\f}$
into a larger family (still defined over $\Q$ albeit on a larger
space) which would contain symmetrizations of all these spaces.
Since the family $\mathscr L_{\l,\f}$ is already conformally
complete (e.g., contains together with each space its conformal
transforms), it is sufficient to symmetrize only with respect to a
single axis, the real line $\R$.

The reflected matrix function $\boldsymbol X^\dag(\mu,\psi,
t)=\overline{\boldsymbol X(\bar\mu,\bar\psi, \bar t)}$ is a
holomorphic matrix function which satisfies the rational
integrable Pfaffian system $\d \boldsymbol X^\dag=\boldsymbol
\Omega^\dag \boldsymbol X^\dag$ on $\P^m\times\P^3\times\P^1$.
This reflection does not affect neither degree, nor the
dimension or the size of the Pfaffian system (obviously keeping
it defined over $\Q$). The block-diagonal matrix function
$$
 \^X(\l,\f,\mu,\psi,t)=\diag\{\boldsymbol X(\l,\f,t),\boldsymbol
 X^\dag(\mu,\psi,t)\}
$$
satisfies the integrable Pfaffian system of
dimension $2\ell$ on the product space birationally equivalent to
$\P^{2m+6}\times \P^1$ with the coordinates $(\l,\f,\mu,\psi,t)$,
\begin{equation}\label{diag-sys}
    \mathrm d\^X=\^\Omega \^X,
    \qquad
    \^X=\begin{pmatrix}\boldsymbol X&\\&\boldsymbol X^\dag\end{pmatrix},
    \quad
    \^\Omega=
    \begin{pmatrix}\boldsymbol \Omega&\\&\boldsymbol \Omega^\dag\end{pmatrix}.
\end{equation}
The corresponding family of subspaces $\^{\mathscr
L}_{\l,\f,\mu,\psi}$ contains all sums $\mathscr L_{\l,\f}+\mathscr
L_{\mu,\psi}^\dag$, in particular, all symmetrizations of the
conformal transforms $(\f^*\mathscr L_{\l})^\ominus_\gamma$.

By its explicit construction, the family \eqref{diag-sys}
(considered as a Pfaffian system on $\P^{2m+6}\times\P^1$) is
integrable, rational and regular. Applied to a system defined over
$\Q$ on $\P^m\times\P^1$, the construction results in a system again
defined over $\Q$ on the larger subspace. The dimension is increased
by the factor of $2$ from $\ell$ to $2\ell$, and the size
(complexity) is increased by an explicit constant factor $O(1)$. To
check that the above doubling preserves quasiunipotence, note that a
small loop (in the sense of Definition~\ref{def:small-loop}) in the
product space $\P^{2m+6}\times\P^1$ projects as a small loop on each
of the components $\P^{m+3}\times\P^1$; the corresponding monodromy
is block diagonal with quasiunipotent (or identical) blocks.

Thus for any Pfaffian system \eqref{pfs} with the matrix
$\Omega=\{\Omega_\l\}$ one can effectively construct its embedding
(as a family) into a larger family with the Pfaffian matrix form
$\^\Omega=\{\^\Omega_{\eta}\}$, $\eta\in\P^{2m+6}$ with the
following characteristic property.

\begin{Lem}\label{lem:conf-completion}
For any parameter value $\l\in\P^m$ there exists a parameter value
$\eta\in\P^{2m+6}$ such that the invariant slope of the derived
operator $D_\l$ is equal to the affine slope of the corresponding
derived operator $\^D_\eta$,
\begin{equation*}
    \forall\l\in\P^m\quad \exists\eta\in\P^{2m+6}\qquad \slope D_\l=\angle\^D_\eta.
\end{equation*}

Thus a uniform bound for the affine slope $\angle\^D_\eta$ of the
family of operators $\^D_\eta$ is at the same time the uniform upper
bound for the invariant slope $\slope D_\l$ of the operators $D_\l$
derived from the initial family.\manualqed
\end{Lem}

\subsection{Proof of the Principal Lemma~\ref{plem:slope}}
\subsubsection{Proving the qualitative part.}\label{sec:firsthalfoftheproof}
To prove the existential finiteness of the latter (Part~\textbf{A}
of the Principal Lemma~\ref{plem:slope}), we need to consider
together with the initial family $\{\Omega_\l\}$, $\l\in\P^m$, and
the corresponding family $D_\l$ of derived operators the conformal
completion $\{\^\Omega_\eta\}$, $\eta\in\P^{2m+6}$, and the
respective family $\{\^D_\eta\}$ as described in
\secref{sec:conf-compl}.

Lemma~\ref{lem:exist-bd} applied to the family $\{\^D_\eta\}$,
guarantees that the affine slope of these operators is bounded
uniformly over $\eta\in\P^{2m+6}$. By
Lemma~\ref{lem:conf-completion}, this means that the invariant slope
of the operators $D_\l$ is bounded by the same constant.

\subsubsection{From qualitative to quantitative bounds.}
It remains to prove Part~\textbf{B} of the Principal Lemma and show
that if the regular family \eqref{pfs} is defined over $\Q$, then
the bound for the invariant slope can be made explicit as follows,
\begin{equation}\label{slope-effect}
    \forall \lambda\in\P_m\ssm\scS
    \qquad\slope D_\l\le s^{(d\ell)^{O(m)}},\quad s=\size{D},\ d=\deg D.
\end{equation}

Indeed, in this case the derived equation $\^D_\eta$ is also defined
over $\Q$ and its size is explicitly bounded by the assertion (3) of
Lemma~\ref{lem:D-construction} (replacing $s$ by its finite power
$s^{O(1)}$ and $\ell$ by $2\ell$ because of the difference between
the families $D$ and $\^D$ does not affect the asymptotic). This
means that the subgraph of the affine slope function
$\sigma(\bullet)$ is a semialgebraic set defined over $\Q$, exactly
as in Example~\ref{ex:K-bound}. Since the slope
$\sigma(\eta)=\angle\^D_\eta$ is bounded by Lemma~\ref{lem:exist-bd}
(cf.~with \secref{sec:firsthalfoftheproof}), Theorem~\ref{thm:basu}
gives then the double exponential upper bound of the form
$s^{(O(1)\cdot d\ell)^{O(2m+6)}}=s^{(d\ell)^{O(m)}}$ for the slope
$\sigma(\bullet)$ on $\P^{2m+6}$. By
Lemma~\ref{lem:conf-completion}, this gives the explicit uniform
bound for the invariant slope $\angle D_\l$ and thus completes the
proof of Principal Lemma~\ref{plem:slope}.

\begin{small}
\begin{Rem}[Proof of Proposition~\ref{prop:inv-slope}]
 \label{rem:inv-slope-finiteness}
Let $D$ be an arbitrary Fuchsian operator. Its conformal transforms
$\f^*D$, and their symmetrizations $(\f^*D)^\ominus_\gamma$ with
respect to all arcs $\gamma=\psi(\R)\conf\R$,
$\f,\psi\in\Aut(\P)\subset\P^3$, constitute a parametric family
$D_\l$, $\l\in\P^3\times\P^3$, which is a regular family of
operators with compact base $\simeq\P^6$ (not necessarily defined
over $\Q$). Application of Lemma~\ref{lem:exist-bd} to this family
proves that the invariant slope $\slope D$ is always finite; this
gives a (very indirect) proof of Proposition~\ref{prop:inv-slope}.

It would be interesting to achieve a direct proof together with an
explicit bound on the invariant slope, e.g., in terms of the order
of the equation, the number of Fuchsian singularities and the
absolute value of the respective characteristic exponents
\cite[Example 19.21]{iy:lade}.
\end{Rem}
\par\end{small}

\section{Counting zeros of functions defined by Fuchsian
equations}\label{sec:renorm}

A linear ordinary differential equation of bounded affine slope
admits an explicit upper bound for the variation of argument of
its nontrivial solutions along paths of bounded length,
sufficiently distant from the singular points of the equation.
For Fuchsian equations, because of the finiteness of the
invariant slope (Proposition~\ref{prop:inv-slope}), such paths
can be drawn with very few restrictions. One can then try and
slit the projective line $\P^1$ along suitably chosen paths
into ``polygonal'' simply connected domains, to which the
argument principle can be applied. A similar approach, also
based on the idea of suitable clustering, was suggested in
\cite{quasialg}, yet its implementation there was conditioned
on the spectral condition imposed on the monodromy group of the
equation. Unfortunately, this condition is algebraically
unverifiable (one cannot, in general, algebraically compute the
monodromy of a regular system along an arbitrary closed loop).
Moreover, this condition in the form required in
\cite{quasialg} does not hold for the Picard--Fuchs system in
general: there are some ``large loops'' whose monodromy does
not possess the necessary spectral properties.

In this section we suggest a way to circumvent this obstacle for isomonodromic families, and
establish explicit upper bounds for the number of zeros under a
weaker condition: the monodromy is required to be quasiunipotent
only around \emph{small loops} (cf.
Definition~\ref{def:small-loop}).

\subsection{Normalized length}\label{sec:norm-length}
We start by introducing some metric characteristics of finite
configurations of $l$ points on the plane $\C$, which are invariant
under the action of the affine group.

Let $T$ be a fixed finite point subset in $\C$.

\begin{Def}\label{def:norm-length}
The \emph{normalized length} of a closed circular arc
$\gamma\subset\C\ssm T$ relative to the finite point set
$T\subset\C$ is the finite positive number
\begin{equation}\label{norm-l}
    |\gamma/T|=\frac1{2\pi}\cdot \frac{|\gamma|}{\dist(\gamma,T)}
\end{equation}
where $|\cdot|$ denotes the usual Euclidean length in $\C$ and
$\dist(\cdot,\cdot)$ the Euclidean distance.

The \emph{normalized length} of a line segment $\gamma$ disjoint
with $T$ is defined as the similar ratio
\begin{equation*}
    |\gamma/T|=\frac{|\gamma|}{\dist(\gamma,T)},
\end{equation*}
differing only by the numeric factor $2\pi\approx6.283\dots$.
\end{Def}

The normalized length clearly is an affine invariant: for any affine
automorphism $\f\in\Aut(\C)$, we have
$|\gamma/T|=|\f(\gamma)/\f(T)|$.

\begin{Ex}
The normalized length of all sufficiently small circles
$\{|t-a|=\e\}$, $a\in T$, $0<\e\ll 1$, is $1$.
\end{Ex}

The normalized length is a crude substitute for the length of a
segment in the hyperbolic (Poincar\'e) metric on the unit disk
considered as the universal covering space for the multiply
connected domain $\C\ssm T$. However, this substitute will be
more convenient for our purposes than the genuine hyperbolic
length, among other things because of the semialgebraicity of
the former.

\subsection{Affine invariant bounds for zeros of solutions
 of differential equations}\label{sec:inv-count}
The notion of normalized length allows to place bounds on the
variation of argument and the number of isolated zeros for
solutions of homogeneous ordinary differential equations in
\emph{affine invariant terms}.

The following statement gives an affine invariant generalization
of Lem\-ma~\ref{lem:var-arg-aff} in terms of the normalized length.

\begin{Lem}\label{lem:var-arg}
Let $D$ be a Fuchsian differential operator of order $k$ with
coefficients of degree $\le d$ and  the singular locus $\S$, and
$\gamma$ a closed circular arc or line segment disjoint with $\S$.

Then the variation of argument of any nonzero solution of the
homogeneous equation $Dy=0$ along the arc $\gamma$ is explicitly
bounded,
\begin{equation}\label{var-arg-bound}
    \Var\Arg y(t)|_\gamma\le k L^{O(d)}\cdot \slope D,
\end{equation}
where $L=|\gamma/\S|$ is the {normalized length} of the arc (see
\eqref{norm-l}), and $\slope D$ the invariant slope of the operator
$D$.
\end{Lem}

\begin{proof}
One can always apply an affine transformation of the affine plane so
that the distance between $\gamma$ and the singular locus of the
operator will be exactly $1$ and the curve itself passes through the
origin. Then the Euclidean length of the curve will not exceed $L$
by the definition of the normalized length, and it will belong to
the circle of radius $L$ by the triangle inequality. The affine
slope $S$ of the operator $D$ in the new chart still is no greater
than $\slope D$ by definition of the invariant slope. Applying
\eqref{var-arg-bd} with $S=\slope D$, $R=|\gamma|=L$ and $r=1$, we
obtain the inequality \eqref{var-arg-bound}.
\end{proof}

In the same way as in the proof of Lemma~\ref{lem:var-arg}, we may
generalize the ``Petrov trick'' to count zeros of solutions of
differential operators in annuli with quasiunipotent monodromy in
affine invariant terms. The starting point is the following Lemma
which is borrowed from \cite{polyfuchs}.

\begin{Lem}[{\normalfont Lemma~10 from
\cite{polyfuchs}}]\label{lem:polyfuchs-annulus}
 If the monodromy of a
\emph{real} differential operator $D$ along the equator of a
symmetric annulus $K=\{\rho_-< |t|< \rho_+\}$ has all eigenvalues on
the unit circle, then the number of zeros of any solution in $K$ is
explicitly bounded,
\begin{equation}\label{pre-annulus-bd}
    \cN(D,K)\le (2k+1)(2B+1),
\end{equation}
where $k$ is the order of the operator and $B=B(D,K)$ the upper
bound for the variation of arguments of any solution of $Dy=0$ along
the boundary circles $C_\pm=\{|t|=\rho_\pm\}$.\manualqed
\end{Lem}

Together with Lemma~\ref{lem:var-arg} above, this estimate
proves the following explicit bound for the number of zeros of
solutions in annuli.

\begin{Lem}\label{lem:annulus-count}
Let $D$ be a Fuchsian operator of order $k$, degree $\le d$ and the
singular locus $\S$, and $K\subseteq\C\ssm\S$ a topological annulus
bounded by two disjoint circles $C_\pm$ (one of which may degenerate to a point).

If the monodromy of $D$ along the equator of $K$ is
quasiunipotent, then the number of isolated zeros of any
solution of $Dy=0$ in $K$ is explicitly bounded\textup:
\begin{equation}\label{annulus-bd}
    \cN(D,K)\le k^2 L^{O(d)}\cdot \slope D,
\end{equation}
where $L$ is the normalized length of the boundary,
$L=|C_-/\S|+|C_+/\S|$.
\end{Lem}

\begin{proof}
We can always make a conformal automorphism which transforms the
annulus $K$ into the annulus bounded by two circles centered at the
origin. Without loss of generality, replacing if necessary $D$ by
its symmetrization around the real axis $\R$, we may assume that $D$
is real. The variation of argument of any solution of the equation
$Dy=0$ along the boundary circles is bounded by $k
|C_\pm/\S|^{O(d)}\slope D\le k L^{O(d)}\slope D$ by
Lemma~\ref{lem:var-arg}. The inequality \eqref{pre-annulus-bd}
then implies the bound \eqref{annulus-bd}.
\end{proof}

\subsection{Admissible systems of slits}
In this section we describe systems of arcs such that slitting
the plane along these arcs subdivides it into components
allowing for application of the counting tools
(Lemmas~\ref{lem:var-arg} and~\ref{lem:annulus-count} from
\secref{sec:inv-count}).

\begin{Def}
The normalized length of a \emph{union} of circular arcs and line segments
$S=\gamma_1\sqcup\cdots\sqcup\gamma_k$, $S\subset\C\ssm T$, is by
definition the sum of the normalized length of all components,
\begin{equation*}
    |S/T|=|\gamma_1/T|+\cdots+|\gamma_k/T|.
\end{equation*}
\end{Def}

\begin{Rem}
Note that for each term $\gamma_i$ above, the normalized length
involves the distance from the set $T$ to $\gamma_i$ and not to
their union $S$. Thus the normalized length \emph{depends} on the
way the set $S$ is represented as a finite union of arcs and segments. In
our constructions, however, this representation will always be clear
from the context.
\end{Rem}

\begin{small}
\begin{Ex}
Let $T$ be any two-point set and $S$ the union of two equal circles
$\gamma_{1,2}$ centered at these points and the shortest line
segment $\gamma_0$ connecting these circles.

If the radii of these circles are equal to $1/3$ of the
distance between the points of $T$, then the normalized length
$|S/T|$ is equal to $3$. This length can be further reduced to
almost $2$ if the radii tend to half the distance between the
points. On the contrary, the normalized length $S/T$ tends to
infinity if the radii of the circles tend to $0$: in this case
$|\gamma_0/T|$ tends to infinity.
\end{Ex}
\end{small}

Let, as before, $T\subset\C$ be a finite point set, and
$S=S_{a,r}\subset\C$ a finite union of circles of the form
$S_{a,r}=\bigcup_{i}\{|t-a_i|=r_i\}$, $a_i\in\C$, $r_i>0$,
$i=1,\dots,k$.

\begin{Def}
The union of circles $S$ is called a \emph{clustering} of the finite
point set $T$, if all these circles are disjoint with $T$ and
pairwise disjoint with each other.
\end{Def}

A clustering subdivides points from $T$ into \emph{nested subsets},
some (or most) of which in principle may be empty or singletons.

Let $S=\bigsqcup_i S_i\subset\C\ssm T$ be a clustering of $T$.

\begin{Def}\label{def:adm-slit}
A finite union $S'=S\sqcup\gamma_1\sqcup\cdots\sqcup \gamma_k$ of
circles $s_i$ and line segments $\gamma_i$ connecting them so that
the circles and segments have only endpoints in common, is called
\emph{admissible system of slits} for a finite point set $T$, if the
complement $\C\ssm\{S'\cup T\}$ consists only of simply connected
domains (of arbitrary shape) and topological annuli bounded by two
circles  (which may degenerate to a circular disk punctured at a
point from $T$).
\end{Def}

\begin{figure}[h]
\centerline{\includegraphics
[width=0.5\textwidth]{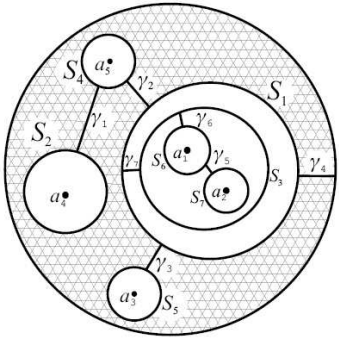}} \caption{Admissible system of slits around a
finite point set.}\label{fig:admissible-slits}
\end{figure}

Clearly, any clustering can be completed to an admissible system of
slits by infinitely many ways. The number of possibilities can be
reduced to finite, if each segment $\gamma_i$ realizes the shortest
path connecting the two respective circles (provided the latter are
not concentric).

\begin{Def}[principal]\label{def:cluster-diameter}
The \emph{cluster diameter} of a finite point set $T\subset\C$
is the infimum of normalized lengths of an admissible system of
slits $S'$ as in Definition~\ref{def:adm-slit}, involving no
more than a given number $c$ of circular arcs:
\begin{equation*}
    \cdiam {T\mid c}=\inf_{S'} \biggl\{\,|S'/T|\:
    S'=\bigsqcup_{i=1}^c S_i\bigsqcup_j \gamma_j  \text{ admissible  for }T\biggr\}.
\end{equation*}
\end{Def}
By this definition, $\cdiam{T\mid c}$ may well be infinite, if the
number of circular slits is too small compared to the number of
points. On the other hand, it is obviously finite if $c$ is
sufficiently large (see the proof of Lemma~\ref{lem:c-semialg}).

\subsection{Admissible system of slits for differential
 equations}\label{sec:slits-oper}
Let $L=p_0(t)\partial^k+\cdots+p_{k-1}(t)\partial
+p_k(t)\in\C[\partial, t]$ be a differential operator with
polynomial coefficients $p_i\in\C[t]$. Its singular locus
$\S_L=\{p_0=0\}$ is a finite point set which will be denoted by $T$.

If $S'$ is an admissible system of slits for the point set $T=\S_L$
in the sense of Definition~\ref{def:adm-slit}, then for any
topological annulus $A\subset\C\ssm S'$ the monodromy operator
$M=M_A$ associated with the equatorial loop (the positively oriented
loop in $A$ which generates $\pi_1(A)$) is defined uniquely modulo
conjugacy $M\mapsto C^{-1}MC$, $\det C\ne 0$.

\begin{Def}
A system of slits $S'$ is \emph{admissible for the operator} $L$, if
it is admissible for the \emph{singular locus} $T=\S_L$ and in addition
the monodromy of $L$ along the equatorial loop of each annulus $A$
appearing in $\C\ssm S'$ is quasiunipotent (cf.~with the assumptions
of Lemma~\ref{lem:annulus-count}).
\end{Def}

The admissible system of slits for a linear system of Pfaffian
equations $\d x=\Omega_\l x$ on a projective line, is defined
analogously. Obviously, the additional constraint imposed by the
requirement on the monodromy map is determined by the solutions, so
that a system of slits is admissible at the same time both for a
regular integrable rational system $\Omega_\l$ on $\P^1$ as in
\eqref{pfs-fib} and for the corresponding derived operator $D_\l$ as
in \eqref{stf}.

As before, we define the \emph{cluster diameter of the singular
locus} of the operator $L$ as the infimum over all admissible
systems of slits involving no more than $c$ circles,
\begin{equation*}
    \cdiam{L\mid c}=\inf_{S'} \biggl\{\,|S'/\S_L|\:
    S'=\bigsqcup_{i=1}^c S_i\bigsqcup_j \gamma_j  \text{ admissible  for
    }L\biggr\}.
\end{equation*}
By construction, $\cdiam {L\mid c}\ge \cdiam{\S_L\mid c}$, since not
all slits admissible for the point set $T=\S_L$ are necessarily
admissible also for $L$: some annuli may have non-quasiunipotent
monodromy.

Let $\Omega$ be an integrable Pfaffian system  on $\P^m\times\P$ and
$D=\{D_\l\}_{\l\in\P^m}$ the associated isomonodromic family of
linear ordinary differential operators in the standard form
\eqref{stf}, equivalent to $\Omega$ in the sense explained in
Lemma~\ref{lem:D-construction}. Denote, as before, by
$p_0\in\C[\l,t]$ the leading coefficient of the family $D$ and by
$\scS=\{\l\:p_0(\l,\bullet)=0\}\subset\P^m$ the degeneracy locus.

Then for any $\l\notin\scS$ we have a uniquely defined finite point
set
\begin{equation}\label{Tl}
    T(\l)=\{t\in\C\:p_0(t,\l)=0\}\subset\C,\qquad \l\notin\scS,
\end{equation}
which consists of at most $l=\deg p_0\le O(\ell^4 d)$ points (some
of which may escape to infinity for particular values of the
parameter $\l$). The linear differential operator $D_\l$ is
nonsingular on $\C\ssm T(\l)$, and we can introduce the function
\begin{equation}\label{Fi}
    \F_c(\l)=\cdiam{D_\l\mid c},\qquad \l\notin\scS.
\end{equation}

\begin{PLem}\label{plem:diam}\indent\par
\textbf{A.} Let $\Omega$ be an integrable rational Pfaffian system
on $\P^m\times\P$ and $D=\{D_\l\}$ the corresponding derived
isomonodromic family of differential operators.

If $\Omega$ is quasiunipotent and $c\ge O(\ell^4 d)$, then the
function $\F_c$ introduced in \eqref{Fi}, is globally bounded
everywhere on $\P^m$.

\textbf{B.} If in addition $\Omega$ is defined over $\Q$ and
$\size\Omega\le s$, then $\F_c$ is defined over $\Q$ and admits an
explicit upper bound,
\begin{equation}\label{Fi-bd}
    \F_c(\l)\le s^{2^{O(d\ell^4m)^5}}\qquad\forall \l\in\P^m\ssm\scS.
\end{equation}
\end{PLem}
In this formulation, as before, $\ell$ is the dimension of the
Pfaffian system $\Omega$, $d$ its degree and $m$ the number of
parameters. The proof of this Lemma occupies sections
\secref{sec:cdiam-alg}--\secref{sec:cdiam-pf}.

\subsection{Semialgebraicity of the cluster diameter}\label{sec:cdiam-alg}
 We start by observing that the cluster diameter of an
isomonodromic family of linear operators is (bounded by) a semialgebraic function
of the parameters. This would be fairly easy to prove using the technique of quantifier elimination if the requirement on the
monodromy was absent in the definition, since the cluster diameter
of a point set is determined by an explicit algorithmic formula. We
show that the isomonodromy is the key to restoring the
semialgebraicity.

\begin{Lem}\label{lem:c-semialg}
In the assumptions of Lemma~\ref{plem:diam}\textbf{A} the function
$\F_c$ is everywhere finite on $\P^m\ssm\scS$ and semialgebraic. The
subgraph of $\F_c$ in $\P^m\times\R_+$ is a semialgebraic set
defined by polynomial equalities and inequalities of degree not
exceeding $(\ell d)^{O(cm)^5}$.

In the assumptions of Lemma~\ref{plem:diam}\textbf{B}, the function
$\F_c$ is defined over $\Q$ and has size  explicitly bounded by a
double exponent, i.e., its graph is defined by real polynomial
(in)equalities with integer coefficients not exceeding
\begin{equation}\label{complexity-Fi}
    \size{\F}\le s^{(\ell d)^{O(cm)^5}}.
\end{equation}
\end{Lem}

\begin{proof}
Denote by $l$ the upper bound for the degree of the coefficients of
the derived equation $D$: by Lemma~\ref{lem:D-construction}, $l$ is
bounded by $O(\ell^4 d)$.

The clusterings of the set $T=\{t_1,\dots,t_l\}$ by $c$ circles
$S_1,\dots,S_c$ can be parameterized by an open subset of the
Euclidean space
\begin{equation}\label{cluspace}
    \scC^{3c}=\C^c\times\R_+^c=\{(a_1,\dots,a_c,r_1,\dots,r_c)\},
\end{equation}
(each circle $S_i$ is defined by the equations $\{|t-a_i|=r_i\}$,
$i=1,\dots,c$).

Consider the product $\P^m\times\scC^{3c}$. The conditions ensuring that
the circles form a clustering of the singular locus
$T(\l)=\S_{D_\l}$ are semialgebraic: the singular points
$t_1,\dots,t_l$, which are algebraic functions of $\l$, should
satisfy the inequalities
\begin{equation*}
    |t_i-a_j|\ne r_j,\qquad
    |a_j-a_k|>r_j+r_k\quad\text{or}\quad|a_j-a_k|<|r_j-r_k|
\end{equation*}
for all roots $t_i$ of the leading coefficient $p_0(t,\l)$ and
all pairs $j\ne k$, $j,k=1,\dots,c$. The latter conditions mean
that the circles are disjoint and non-nested (resp., disjoint
and nested). We add to these conditions the inequalities
$a_i\ne a_j$ which will guarantee that the shortest slits
connecting any two given circles, are uniquely defined (as no
two circles are concentric), and the conditions that the ratios
$(t_i-a_j)/(t_i-a_k)$ are non-real (this will guarantee that
the shortest slits will not pass through the singular point
$t_i$).

The points $t_i$ themselves are defined by the equation
$p_0(\l,t_i)=0$ (roots of the leading coefficient of the
differential operator $D_\l$) and the degree of that leading
coefficient should be maximal (equal to $l$) to avoid escape of the
roots at infinity.

Altogether we have $O(c+l)^3$ real algebraic equalities/inequalities
of degree $\le l$ and size $\le s$ in the space of $3c+2m=O(c+m)$
real variables $\P^m\times\scC^{3c}$.

Admissible clusterings (defined by collections of circles, without
any reference to the monodromy yet) form a semialgebraic open subset
of the total space $\P^m\times\scC^{3c}$, which consists of a large
number of connected components (different topological configurations
of singular points and circles). Each component $C_\alpha$ is a
semialgebraic set by \cite[Theorem 16.13]{roy.e.a.}. Moreover, its
algebraic complexity can be explicitly bounded: each connected
component can be defined over $\Q$ by polynomial (in)equalities of
degree not exceeding $l^{O(c+m)^3}$ and size (complexity) at most
$s^{l^{O(c+m)^3}}$.

The admissible system of slits $S'$ based on a given clustering
$S\in\scC^{3c}$, can be achieved by a finite number of choices
(the number does not exceed $O(c^2)$) of the segments
$\gamma_j$. Indeed, one can always assume that the segments
$\gamma_j$ realize the shortest distance between each pair of
circles of the clustering (see
Fig.~\ref{fig:admissible-slits}), and by construction this
choice is uniquely defined over each connected component
$C_\alpha\subset\P^m\times\scC^{3c}$. For any choice the
normalized length of the resulting system will be a
semialgebraic function on the position of the singularities and
the clustering parameters (the coordinates of the corresponding
point in $\scC^{3c}$). In other words, the normalized length of
any system of slits $|S'/T(\l)|$ can be considered as a
multivalued semialgebraic function
\begin{multline}\label{Fi2}
    \F_c(\l,S)=\{|S'/T(\l)|\:S'=S\sqcup\gamma_1\sqcup\cdots\sqcup\gamma_{2l}\},
    \\
    \quad\gamma_j\text{ shortest slits},
\end{multline}
on $\P^m\times\scC$. Each of the finitely many continuous
branches of this function is semialgebraic on each connected
component $C_\alpha\subset\P^m\times\scC$.

Moreover, each continuous branch of this function can be
majorized by a function defined over $\Q$, if the initial
system is defined over $\Q$. Indeed, the distance between a
point $t_i$ and the circle $S_j=\{|t-a_i|=r_i\}$ is
$\bigl|r_i-|t-a_i|\bigr|$, the distance between two circles is
given by a similar formula, all of them involving only the
coordinate variables, the absolute value and the coefficients
$0,\pm 1$. Thus the normalized length of all segments is
defined over $\Q$. In the same manner the normalized length of
the circles is defined over $\Q$ (this is the reason why the
factor $2\pi$ appears in Definition~\ref{def:norm-length} of
the normalized length). The complexity of the formula defining
$\F_c(\l,S)$ is at most polynomial in $m+c$, since all
coefficients are bounded by a common constant $O(1)$ (for
brevity we denote the majorant by the same symbol as the
initial function $\F_c$).

It remains to take into account the requirement on the
monodromy of the system (equation). In general, the
\emph{monodromy of a linear equation depends in a non-algebraic
way on the equation}, thus the admissibility of a system of
slits \emph{cannot be defined by an algebraic condition}.
However, since the initial system is integrable, \emph{the
monodromy is constant along any continuous branch
\textup(system of slits\textup) based on the clustering varying
over a connected component} $C_\alpha$. In other words, the
requirement of quasiunipotence reduces to a \emph{branch
selection of the function $\F_c(\l,S)$, defined on
$\P^m\times\scC$}: some of the branches give the normalized
length of an admissible system of slits, while others do not.

It may well happen that for a given clustering $S$ there is no
admissible system of slits based on this clustering, in which
case we set $\F_c(\l,S)=+\infty$. In the case where several
systems of slits based on the same clustering are admissible,
we can choose any of them to evaluate $\F_c(\l,S)$, or use the
minimal value. This will not affect the complexity of the
function $\F_c(\l,S)$.

Ultimately we can express the (majorant for the) cluster diameter
\eqref{Fi} as an infimum of a semialgebraic function,
\begin{equation}\label{Fi-semialg}
    \F_c(\l)=\inf_S\{\F_c(\l,S)\:(\l,S)\in\P^m\times\scC^{3c}\}\le +\infty,
\end{equation}
which itself is semialgebraic by the Tarski--Seidenberg theorem
(quantifier elimination principle). Moreover, since the complexity
of the quantifier elimination algorithm is known, we can guarantee
that the polynomial (inequalities) defining the graph of $\F_c(\l)$
over $\Q$, have degree at most $l^{O(c+m)^3O(c)O(m)}\le l^{O(cm)^5}$
and size (complexity) at most $s^{l^{O(cm)^5}}$. Substituting the
value $l=O(\ell^4 d)$, we obtain the bound \eqref{complexity-Fi}.

It remains to show that the function $\F_c$ takes finite values
for all values of the parameter $\l\notin\scS$, i.e., for
matrices $\Omega_\l$ (resp., operators $D_\l$) with finite
singular locus, provided that $c\ge l$.

In this case the clustering $S_\l$ which consists of exactly
$l$ circles centered at each singular point $t_i\in T_\l$ and
having sufficiently small radius $r_i\ll 1$, can be completed
by finitely many segments to an admissible system of slits.
Indeed, one has to connect the small disks in an arbitrary way
with each other and with a large circle ``centered at
infinity'' (the circle whose exterior contains only one
singular point at $t=\infty$): the only annuli that are formed
by these slits, are punctured disks around singularities, and
their monodromy is quasiunipotent by the quasiunipotence
assumption on the initial system $\Omega$ (cf.~with
Remark~\ref{rem:small-loops-1d}). Thus $\F_c(\l,S_\l)<+\infty$,
hence $\F_c(\l)<+\infty$.
\end{proof}

Note, however, that the finiteness of values of the function
$\F_c(\l)$ for $\l\notin\scS$ does not imply yet its local
boundedness. This last step is achieved in the next section.

\subsection{Local boundedness of the cluster diameter in
 one-paramet\-ric quasiunipotent families}\label{sec:cdiam-bd}
The cluster diameter of a finite point set $T=T(\l)$ (resp. a
family $D_\l$ of equations) depending on a parameter, remains a
continuous (hence locally bound\-ed) function of $\l$ as long
as the points of $T(\l)$ (resp. the singular points of $D_\l$)
do not collide. In an analytic collision of two or more points
the topological structure of the underlying clusters must be
chosen depending on the relative ``speed'' of the colliding
points. An explicit choice of this structure is possible in
one-parametric families.

\begin{Lem}\label{lem:c-bound}
Let $\Omega=\{\Omega_\l\}_{\l\in(\C^1,0)}$ be an integrable
meromorphic quasiunipotent system on $(\C^1,0)\times\P^1$ and
$D=\{D_\l\}$ the corresponding derived equations with rational
coefficients of degree $\le l$.

Then the function $\F_c(\l)=\cdiam{D_\l\mid c}$ is bounded over all
$\l\in(\C^1,0)$, if $c\ge 3l$.
\end{Lem}

\begin{proof}
Consider the singular locus of the system: in the coordinates
$(\l,t)\in(\C^1,0)\times\C^1\subset(\C^1,0)\times\P^1$ it is given
by an equation $p(\l,t)=0$, polynomial in $t$ with coefficients,
holomorphic on $\l$ of degree $\le l$. Without loss of generality we
may assume that $p(0,\bullet)\ne0$ (otherwise divide $p$ by a
suitable power of $\l$).

The equation has $l$ roots $t_1(\l),\dots,t_l(\l)$, which are
algebraic functions on $\l$ and, as such, can be expressed by
converging Puiseaux series. Passing to a fractional power of the
parameter $\l^{1/d}=\e$, we may assume that each root is a
holomorphic function, $t_j=t_j(\e)$, of the parameter $\e$. The
quasiunipotence of the system is preserved by such
re-parametrization.

We will construct a clustering
$S_\e=\{(a_i(\e),r_i(\e))\}\in\scC$ of all sets
$T(\e)=\{t_1(\e),\dots,t_l(\e)\}$ for all sufficiently small
values of $\e$ and an associated admissible system of slits for
$S_\e$ (a continuous branch of the function $\F_c(\e)$ in the
terminology of \secref{sec:cdiam-alg}) such that the function
$\F_c(\e,S_\e)$ will be finite as $\e\to0$ (see \eqref{Fi2}).

\begin{figure}
\centerline{\includegraphics[width=0.8\hsize]{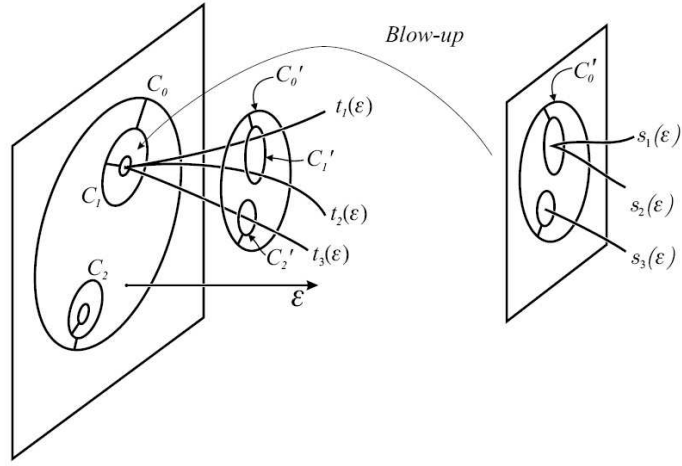}}
\caption{Construction of the admissible system
of slits.}\label{fig:slits}
\end{figure}

1. The outermost circle $C_0$ of the clustering is the circle which
contains all points of the set $T(0)$ and is of distance at
least $1$ from them.

The next embedded level is the union of circles $C_j$ of radius $\rho/2$
centered at all distinct points of the set $T(0)$, where
\begin{equation}
   \rho=\min_{t_i\neq t_j} \{\abs{t_i-t_j}, \quad t_i\in T(0)\}.
\end{equation}
Clearly, all
these circles will be disjoint with $T(\e)$ for all sufficiently
small $\e$, and the normalized length relative to $T(\e)$ of these
circles remains bounded as $\e\to0$. We can add shortest slits
between the outermost circle $C_0$ and some of the first level
circles $C_1$ to make the complement simply connected.

Construction of the next level circles is organized in the same way
relative to circles of the first level, so we will explain it only
for the circle $C_1$ around one of the points $t_1(0)\in T(0)$,
assuming for simplicity that this point is at the origin,
$t_1(0)=0$, so that $C_0=\{|t|=1\}$.

2. If among the roots $t_j(\e)$ there is only one such that
$t_j(0)=0$, i.e., if the origin is a ``simple'' (non-multiple) point
of $T_0$, then the construction in this circle stops, and the
(degenerate) annulus $\{0<|t|<1\}$ has finite relative length for
all small $\e$. The monodromy along this annulus is quasiunipotent
by Remark~\ref{rem:small-loops-1d}.

3. If there is more than one root $t_j$ with $t_j(0)=0$, then
several holomorphic functions $t_j\:(\C^1,0)\to\C$ have the same
$0$-jet. Assume that these functions are labeled as
$t_1(\e),\dots,t_p(\e)$, $p\ge 2$. Let $k\ge 1$ be the first natural
number such that $k$-jets of $t_1,\dots,t_p$ (in $\e$) are \emph{not
all equal} between themselves.

After the rescaling $t\mapsto s=(t-t_1(\e))/\e^k$ in the new local
chart $s$ we will obtain $p$ functions $s_1(\e),\dots,s_p(\e)$,
still holomorphic in $\e$, but with the limits $s_j(0)$ not all
coinciding.

Construct a circle $C'_0$ which in the chart $(\e,s)$ is large
enough to encircle all points $s_j(0)$ and has distance at least $1$
from them, and the smaller circles $C_j'$ of the form
$\{|s-s_j(0)|=\tfrac13\rho'\}$ centered at each distinct point of
the set $T'=\{s_1(0),\dots,s_p(0)\}$, where
\begin{equation}
   \rho'=\min_{s_i\neq s_j} \{\abs{s_i-s_j}, \quad s_i\in T'\}.
\end{equation}
In the original chart $t$
these will be very small circles (of radius $O(\e^k)$). By
construction, the normalized length of $C_0',C_1',\dots,C_p'$
depends only on the position of the points inside $C_1$, since all
other points of $T(\e)$ are incomparably far. On the other hand, because
the normalized length is affine invariant, it can be computed in the
chart $s$, in which it is finite uniformly over all $\e\to0$ by the
same arguments as in the step 1 of the proof.

The circles $C_0',C_1',\dots C_p'$ will be included in the
clustering; the admissible system of slits is complemented by the
slits between $C_0'$ and some of the $C_j'$ to make the slit
interior of $C_0'\ssm\bigcup_{j=1}^p D_j'$ simply connected (here
$D_j'$ are the disks bounded by the circles $C_j'$). As in step
1, these extra slits will have uniformly finite normalized length as
$\e\to0$.

%
%

4. We need to show that the annulus bounded by $C_1$ and $C_0'$ has
an admissible (quasiunipotent) monodromy. This is evident if the
fiber $\{\e=0\}$ is not in the singular locus of $\Omega$, since
then this monodromy coincides with the monodromy of $C_1$, which is
admissible by the Kashiwara theorem~\ref{thm:kashiwara}. In the
opposite case one needs slightly more involved arguments.

\begin{figure}
\centerline{\includegraphics[width=0.97\hsize]{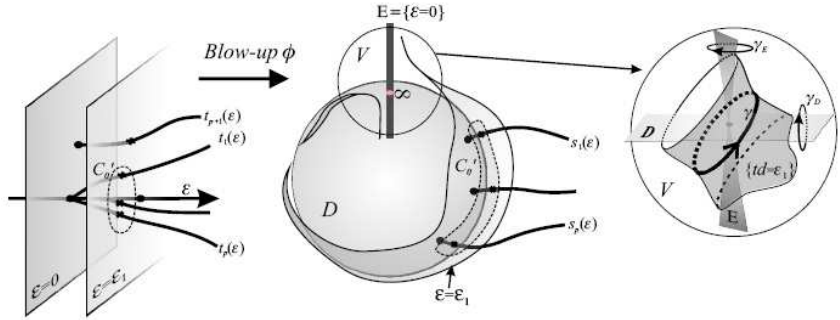}}
\caption{Quasiunipotence of the loops encircling a cluster.}\label{fig:QUP}
\end{figure}

\begin{small}
\begin{Lem}
The annulus bounded by $C_1$ and $C_0'$ on $\{\e=\e_1\}$, with
sufficiently small $\e_1$, has an admissible (quasiunipotent)
monodromy.
\end{Lem}
Note that this is not an immediate corollary of the Kashiwara
theorem, as the circle $C_0'$ has zero linking number with
$\{\e=0\}$, so it cannot be boundary of a holomorphic disc with
center at the origin.

\begin{proof}
First, applying the translation $t\to t-t_1(\e)$, we can assume that
one of our curves coincides with the axis $\{t=0\}$.

Consider first the case where the number $k$ which appeared on Step
3 above is equal to $1$: this means that among the singularities
forming the cluster, there are at least two points strictly
$O(\e)$-distant from each other as $\e\to0$,
$|t_j(\e)-t_i(\e)|^{-1}=O(\e^{-1})$.

Consider the blow-up $\phi$ given in the affine chart $\e\ne0$ by
$(t,\e)\to (\e, s=t/\e)$, and denote by $D\simeq\P^1$ the
corresponding exceptional divisor, see Fig.~\ref{fig:QUP}. The
lifting of each curve $t=t_j(\e)$ is the curve $s=s_j(\e)$, with
curves corresponding to our cluster (i.e. for $j=1,\dots,p$) tending
to some well-defined limit on $D$ as $\e\to 0$, and other curves not
intersecting some neighborhood $U$ of $D$. The lifting of $C_0'$ is
a circle on $\{\e=\e_1\}$ (still denoted by $C'_0$) encompassing all
points $(s_j(\e_1),\e_1)$.

Now, $C_0'$ is homotopic along the leaf $\e=\e_1$ and inside $U$  to
a small circle around $s=\infty$. Consider the second affine chart
of the blow-up covering the neighborhood of $\e=0$. In this affine
chart the blow-up is given by the formulas $(t,\e)\to(t, d=\e/t)$,
the leaf $\{\e=\e_1\}$ is defined by the equation $\{td=\e_1\}$, and
$C_0'$ is homotopic to the curve $\gamma=\{\left(t=\exp(2\pi i
\theta), \d= \e_1\exp(-2\pi i \theta)\right),\ \theta\in[0,1]\}$.
Again, this curve has linking numbers of different signs with the
$t$-axis and $d$-axis, so it cannot be a  boundary of a holomorphic
disc passing through the point $(t,d)=(0,0)$, i.e. it is not a small
loop. In fact, if we denote by $\gamma_D=\{(t=\exp(2\pi i \theta),
d=d_0)\}$ and $\gamma_E=\{\left(t=t_0, d= \e_1^k\exp(2\pi i
\theta)\right)$ the two small loops encircling $D=\{t=0\}$ and
$E=\{d=0\}$, then $\gamma=\gamma_D \gamma_E^{-1}$ in
$\pi_1(V\ssm\{td=0\})$, where $V$ is a small neighborhood of the
point $(t,d)=(0,0)$.

Consider the lifting $\phi^*\Omega$ of  the connection $\Omega$.
This lifting has admissible monodromy along small loops around the
strict transform of the singular locus of $\Omega$. It also has
admissible monodromy along small loops around $D$, since their
projections are still small loops, so their monodromy is still
admissible by Kashiwara theorem. This in means particular that the
monodromies $M_D$ and $M_E$ of $\phi^*\Omega$ along $\gamma_D$ and
$\gamma_E$ respectively, are quasiunipotent. But $M_D$ and $M_E$
commute since $D$ and $E$ form a normal crossing, so that
$\pi_1(V\setminus\{td=0\})$ is commutative. Thus the monodromy along
$\gamma$ is equal to the product $M_D M_E^{-1}$, which is
quasiunipotent as asserted.

For $k>1$ one should perform $k$ blow-ups in order to get the same
situation near the last exceptional divisor. Again, the only curves
whose strict transforms will intersect the last divisor $D$ will be
the curves corresponding to our cluster, and, deforming the loop
encircling them to a neighborhood $V$ of the point of intersection
of $D$ with the previous exceptional divisor, we represent it as a
product of two monodromies. Both monodromies are quasiunipotent due
to the Kashiwara theorem, and they commute for topological reasons,
so their product is again quasiunipotent.
\end{proof}
\end{small}

5. One can further iterate this construction, applying it to
multiple points of the set $T'(0)$, if any, and constructing circles
of the second level of embedding. Clearly, the maximal multiplicity
goes down by at least one in each descent step, thus the tree-like
clustering process terminates (in each branch) no later than after
$l$ steps, the total number of the circles in the clustering being
at most $3l$. Thus for $c\ge 3l$ the conformal diameter
$\F_c(\l)$ is bounded over all $\l\in(\C^1,0)$ as claimed.

\end{proof}

\begin{Rem}
The construction of the admissible system of slits for a
quasiunipotent integrable system is completely classical. The reader
will easily recognize in it the ``screens'' of Fulton and MacPherson
\cite{fulton-macpherson} and/or the desingularization algorithm used
in the proof of the Kashiwara theorem \cite{kashiwara:qu}.
Definition~\ref{def:cluster-diameter} of the cluster diameter which
is affine invariant was constructed in light of these two proofs.
\end{Rem}

\subsection{Demonstration of the Principal Lemma~\ref{plem:diam}}\label{sec:cdiam-pf}
Consider an integrable rational family of quasiunipotent systems
$\Omega=\{\Omega_\l\}$ on $\P^m\times\P^1$ with the deleted
indeterminacy locus $\scS\subset\P$.

The cluster diameter $\F_{c}(\l)$ is a semialgebraic function
on $\P^m\ssm\scS$, continuous (hence locally bounded) outside
$\scS$, by Lemma~\ref{lem:c-semialg}.

We show that this function is locally bounded at each point of
$\P^m$ (including points of $\scS$) for $c\ge3l$. Indeed, assuming that $\F_c$ is
unbounded, by the Curve Selection Lemma we can assume that $\F_c$
grows to infinity along some real analytic (even algebraic) curve in
$\P^m\ssm\scS$ with an endpoint in $\scS$ (cf.~with the proof of
Lemma~\ref{lem:exist-bd}). However, this is impossible by virtue of
Lemma~\ref{lem:c-bound}.

Thus $\F$ is locally bounded near each point of $\P^m$. Because of
the compactness of the latter, we conclude that the cluster diameter
is globally bounded, as asserted in Part \textbf{A} of the Lemma.

To prove Part \textbf{B}, note that by the second assertion of
Lemma~\ref{lem:c-semialg}, the function $\F_c$ is defined over
$\Q$ by polynomial (in)equalities of degree  not exceeding
$(\ell d)^{O(lm)^5}$ and its complexity is explicitly bounded
by the double exponent $s^{(\ell d)^{O(lm)^5}}$ as in
\eqref{complexity-Fi}, where $l=O(\ell^4 d)$. By the
``quantization principle'' (Theorem~\ref{thm:basu}, cf.~with
Example~\ref{ex:K-bound}), the uniform maximum of $\F_c$ does
not exceed the double exponential expression
\begin{equation*}
    \biggl(s^{(\ell d)^{O(lm)^5}}\biggr)^{(\ell d)^{O(lm)^5}}
    =s^{(\ell d)^{O^+(d\ell^4m)^5}}
    =s^{2^{O^+(d\ell^4m)^5}},
\end{equation*}
which coincides with \eqref{Fi-bd}.\qed

\subsection{Proof of Theorems~\ref{thm:exist} and~\ref{thm:count}}
Consider the Pfaffian system \eqref{pfs} on $\P^m\times\P^1$ and the
corresponding derived family $D$ of homogeneous differential
equations \eqref{famlodo}, and fix an arbitrary value of the
parameters $\l\notin\scS$.

The invariant slope of the corresponding operator $D_\l$ is
explicitly bounded by the Principal Lemma~\ref{plem:slope}, by a
double exponential expression \eqref{slope-explicit}:
\begin{equation}\label{bound-D}
    \slope D_\l\le s^{(d\ell)^{O(m)}}.
\end{equation}
In addition, the cluster diameter of the family $D_\l$ with
$c=3l=O(\ell^4 d)$ circular slits is uniformly bounded, and
does not exceed the double exponential bound
\eqref{Fi-bd}
\begin{equation}\label{bound-L}
    L=\cdiam{D_\l\mid c}\le s^{2^{O^+(\ell^4dm)^5}},
\end{equation}
by Principal Lemma~\ref{plem:diam}.

For each of the simply connected domains  $U_i$ formed by the
admissible slits, the variation of argument does not exceed the
product \eqref{var-arg-bound} (Lemma~\ref{lem:var-arg}) which,
after substitution of the bounds \eqref{bound-D}
and~\eqref{bound-L} and $k=O(\ell^2)$ yields the double
exponential bound $s^{2^{O^+(d\ell^4m)^5}}$ which absorbs all
other terms in the product. This places an upper bound for the
number of zeros $\cN(\Omega/U_i)$. By the same token, the same
bound holds for the number of zeros $\cN(\Omega/U_j)$ in each
annulus $U_j$ formed by the admissible slits: this follows from
Lemma~\ref{lem:annulus-count}.

Since any triangle $T$ may intersect at most  $O(l)$ different
domains (recall that $l$ measures the topological complexity of
the singular locus and hence the number of simply connected
domains/annuli, appearing by admissible slits) the total number
of zeros $\cN(\Omega)$ is bounded by the double exponent as
asserted.

\subsection{Concluding remarks}\label{sec:concl}
The growth of the upper bound for $L$ (the cluster diameter), which
is the fastest of the three terms in the products
\eqref{var-arg-bound} and \eqref{annulus-bd}, is determined by the
complexity of the clustering algorithm. The asymptotics can be
improved at the cost of transparency.

For instance, among the $l=O(\ell^4 d)$ singular points of the
derived equation, almost all (except for $d$) are \emph{apparent
singularities} at which solutions of the equation are non-ramified
(and even remain holomorphic). Clearly, the location of such apparent
points cannot affect the quasiunipotence of the monodromy along an
annulus. On the other hand, the normalized distance from these
points to the slits of an admissible system is important. Thus
instead of the $3l$ slits used in the proof, one can use only $3d$ slits,
while the corresponding normalized length will be given by a formula
which is considerably more complex (involving distance between roots
of an equation and a given circle). Yet since the complexity plays
much less crucial role than the number of variables, the overall
result will be ultimately better.

Other more subtle modifications can be made in order to prove a bound which is double exponential only in the number of parameters of the original system.
However, for the purpose of preserving the clarity of exposition, we shall
not investigate the necessary modifications in this paper.

Another interesting example is that of hyperelliptic integrals. In
this case the monodromy of \emph{any} annulus free from the singular
points, is quasiunipotent. This follows from the Lyashko--Looijenga
theorem \cite{looijenga} asserting that any deformation of critical
values of a univariate polynomial can be achieved by a suitable
deformation of its coefficients (a fact which is wrong for
multivariate polynomials). Thus any $k< n=\deg H$ singular points of
a hyperelliptic Picard--Fuchs system can be isomonodromically
deformed into a degenerate singularity, so any loop is homotopic to
a small loop.

This observation means that in the construction of the admissible
system of slits one can drop the monodromy condition and hence skip
the step of isolating a connected component $C_\alpha$ in the proof
of Lemma~\ref{lem:c-semialg}. The problem of optimizing the
admissible system of slits becomes the problem from the elementary
Euclidean geometry on the plane, namely the computation of the cluster
diameter of a point set $T$ without additional restrictions. Using
elementary arguments, one can estimate the cluster diameter by a
simple exponent of $l=\#T$: $\cdiam{T\mid 3l}\le 2^{O(l)}$.

Given that the hyperelliptic Hamiltonian of degree $n+1$ has
$n$ coefficients (and not $O(n^2)$, as a bivariate polynomial),
we have a better bound for the invariant slope of the
corresponding derived equation. The corresponding double
exponent $2^{2^{O^+(n)}}$ absorbs all other dependencies and
altogether the mentioned improvements give the bound described
in Theorem~\ref{thm:hyperel}. We leave the details to the
reader.

\appendix
\section{Appendix. Complexity of the Picard--Fuchs system}\label{sec:compl-pf}

In this appendix we show that the Picard--Fuchs system of linear
Pfaffian equations satisfied by the periods of monomial forms is
defined over $\Q$ and has an explicitly bounded complexity (size)
and prove Theorem~\ref{thm:pf-compl}.

This proof can be achieved by inspection of the effective derivation
of the Picard--Fuchs system in \cite{bullscimat-2002}, see also
\cite{redundant} for an earlier version. For the readers'
convenience, we reproduce the construction here together with all
required estimates.

\subsection{Effective decomposition in the Petrov module}
In what follows we fix a natural number $n\in\N$ and denote by $H$
the polynomial
\begin{equation*}
    H(x,\l)=\sum_{0\le |\alpha|\le n+1}\l_\alpha
    x^\alpha\in\Z[x,\l],
    \qquad x=(x_1,x_2),\ \l\in\C^{m+1}
\end{equation*}
(we use the standard multiindex notation, $x=(x_1,x_2)$,
$\alpha=(\alpha_1,\alpha_2)\in\Z_+^2$). For each specification of
$\l\in\C^{m+1}$ we obtain a complex polynomial
$H(\bullet,\l)\in\C[x]$, which for a generic value of $\l$ is a
Morse function with the principal homogeneous part $\^H(\bullet,\l)$
having an isolated critical point of multiplicity $n^2$ at the
origin, where
\begin{equation*}
    \^H(x,\l)=\sum_{|\alpha|=n+1}\l_\alpha x^\alpha
\end{equation*}
(the latter condition occurs if and only if $\^H(\bullet,\l)$ is square-free in
$\C[x]$).

It is well-known that the monomials $x^\alpha$, $0\le
\alpha_{1,2}\le n-1$, constitute a basis in the quotient space
\begin{equation}\label{quot-space}
    Q_\l=\C[x]/I_\l,\qquad I_\l=\biggl<\pd{\^H}{x_1}(x,\l),
    \pd{\^H}{x_2}(x,\l)\biggr>\subset\C[x_1,x_2]
\end{equation}
by the gradient ideal $I_\l$ for almost all (though not all) $\l$
such that $\^H(\bullet,\l)$ is square-free. For such values of the
parameters $\l$ the forms $\mu_\alpha=x^\alpha\,\mathrm d
x_1\land\mathrm dx_2$ form the basis of the module of the Brieskorn
lattice and their (monomial) primitives generate the so called
Petrov module \cite[\parasymbol 26E]{iy:lade}. More precisely, we
have the following result \cite[Theorem~26.21]{iy:lade}.

Let $\omega_\alpha\in\L^1(\C^2)$ be the monomial 1-forms such
that $\d\omega_\alpha=\mu_\alpha$, $0\le\alpha_{1,2}\le n-1$,
and $\omega$ any other monomial form.

\begin{Thm}[{\normalfont see \cite{iy:lade}}]\label{thm:petrov-mod}
If the monomials $x^\alpha$ as before generate the quotient space
\eqref{quot-space} for a given value $\l$, then for any monomial
form $\omega$ there exist univariate polynomials $p_\alpha\in\C[t]$
and bivariate polynomials $u,v\in\C[x_1,x_2]$ such that
\begin{equation}\label{rel-coh}
\begin{gathered}
    \omega=\sum_{\alpha}(p_{\alpha}\circ H)\cdot \omega_\alpha +u \d
    H+\d v\quad 0\le\alpha_{1,2}\le n-1,\\
    p_{\alpha}\in\C[t],\ u,v\in\C[x_1,x_2],\quad
    \d H=\pd H{x_1} \d x_1+\pd H{x_2}\d x_2,
    \\
    \left\{\begin{aligned}
    &(n+1)\deg p_{\alpha}+\deg\omega_\alpha\\
    &\deg v\\
    &n+\deg u
    \end{aligned}
    \right.
    \quad\le\deg\omega.
\end{gathered}
\end{equation}

In a similar way, the forms $\mu_\alpha=\d\omega_\alpha$
themselves generate all polynomial $2$-forms $\L^2$ as a module
over $\C[t]$ modulo the submodule $\d H\land\L^1$\textup: any
monomial 2-form $\mu$ admits a representation
\begin{equation}\label{rel-coh2}
    \mu=\sum_\alpha (p_\alpha\circ H)\cdot\mu_\alpha+\d H\land\eta,\qquad
    \eta\in\L^1,
\end{equation}
with analogous inequalities between the degrees of the coefficients
$p_\alpha$ and the ``incomplete ratio'' $\eta$.
\end{Thm}

While the Theorem says nothing about the dependence of the
result of the division on $\l$, we claim that, as functions of
$\l$, the polynomials $p_{\alpha}$ (the remainders) and the
1-form $\eta$ are rational and defined over $\Q$ if the left
hand sides $\omega$, resp., $\mu$ are over $\Q$. For future
calculations we will need only the situation when the expanded
forms $\omega$, $\mu$ are monomial of degree not exceeding
$O(n^2)$.

\begin{Prop}\label{prop:size-division}
Assume that $\omega$ \textup(resp., $\mu$\textup) is a monomial
1-form \textup(resp., 2-form\textup) of degree at most
$O(n^2)$.

Then one can construct an expansion \eqref{rel-coh}, resp.,
\eqref{rel-coh2}, so that the functions
$p_{\alpha}\in\C[t](\l)$, the polynomials $u,v\in\C[x,y](\l)$
and the polynomial 1-form $\eta\in\C(\l)\otimes\L^1[x,y]$ are
all defined over $\Q$, their degrees are bounded by explicit
expressions growing no faster than $O(n^2)$ and their sizes do
not exceed $2^{O(n^3)}$.
\end{Prop}

\begin{proof}
Both systems \eqref{rel-coh}, \eqref{rel-coh2} are linear with
respect to the unknown polynomials $p_{\alpha}(t)$, $u(x),v(x)$ and
polynomial 1-form $\eta$, thus the latter could all be found using
the method of indeterminate coefficients.

The degrees of these objects are explicitly bounded by
Theorem~\ref{thm:petrov-mod}, so we have an explicit control over
the number of unknown indeterminate coefficients to be found,
\begin{equation*}
    \deg p_\alpha=O(1),\quad \dim\{p_\alpha\}=O(n^2),\quad \dim u,\dim v,\dim\eta =O(n^2),
\end{equation*}
(where by $\dim(\bullet)$ we mean the number of the unknown
coefficients in the expansion of these objects as polynomials in $t$
and $x$ respectively). Altogether we see that each of the systems
\eqref{rel-coh}, \eqref{rel-coh2} reduces to a system of linear
(non-homogeneous) algebraic equations. The number of unknowns $N$ of
this system is at most $O(n^2)$, and all entries of the
corresponding matrix are polynomials from $\Z[\l]$: these
polynomials are computed using the explicit expression for $H$ and $\d H$
respectively. The degrees of the entries in $\l$ are at most $\deg
p_\alpha\le O(1)$, since $H$ is linear in $\l$ by assumption. The
size (complexity) of the entries is also bounded by $n$ (the biggest
natural number that appears in the expansion of $\d H$ as a function
of $\l$).

Solutions of such a linear system can be obtained as in
Example~\ref{ex:linsys}. By construction, the corresponding minors
will be polynomials from $\Z[\l]$ of degrees not exceeding
$O(N)=O(n^2)$ and the size (complexity) not exceeding $N!\,n^{N}\le
2^{O(n^3)}$.
\end{proof}

\begin{Rem}
A more accurate analysis carried out in \cite{bullscimat-2002} shows
that the denominators of the rational fractions representing the
polynomials $p_\alpha$ may involve only the parameters $\l_\alpha$
with $|\alpha|=n+1$ corresponding to the principal homogeneous part
of $H$; the dependence on the non-principal coefficients with
$|\alpha|\le n$ is always polynomial.
\end{Rem}

\subsection{Gelfand-Leray derivative}
For any fixed (independent of $\l$) polynomial 1-form
$\omega\in\L^1(\C^2)$ and any multiindex $\alpha$ the derivative
\begin{equation}\label{int-mono}
    \frac{\partial}{\partial \l_\alpha}\oint_{\delta}\omega,\qquad
    \delta\subset\{H=0\},\quad H=\sum_{0\le |\alpha|\le n+1} \l_\alpha x^\alpha
\end{equation}
can be expressed as the integral of another (in general, only
rational in $x$) 1-form $\eta$ over the same cycle
$\delta\subset\{H=0\}$, if the latter satisfies the identity
\begin{equation}\label{g-l}
    x^\alpha\d\omega=-\d H\land\eta.
\end{equation}
This follows easily \cite[Lemma~3]{bullscimat-2002} from the
``standard'' Gelfand--Leray formula which corresponds to
$\alpha=(0,0)$ \cite[Theorem~26.32]{iy:lade}. We use this
observation to express the derivative of the period matrix $X$ of
the monomial 1-forms $\omega_\alpha$ forming the basis.

\subsection{Effective derivation of the Picard--Fuchs system}
In this section we complete the proof of Theorem~\ref{thm:pf-compl}.

Let $X$ be the period $\ell\times\ell$-matrix of monomial 1-forms
$\omega_\alpha$, $0\le\alpha_{1,2}\le n-1$, $\ell=n^2$, cf.~with
\eqref{X}. Denote as before $\mu_\alpha=\d \omega_\alpha$.

For any $\mu_\alpha$ the multiple $H\cdot\mu_\alpha$ is a polynomial
(in $x$) 2-form of degree $\le (n+1)+n^2$ with coefficients
polynomially depending on $\l$. By \eqref{rel-coh2}, there exist
decompositions (the results of division by $\d H$ with remainder)
\begin{equation}\label{pf-dec-1}
    H\cdot\mu_\alpha=\sum_{0\le |\beta|\le n+1}
    (P^\star_{\alpha\beta}\circ H)\cdot \mu_\beta+\d H\land\eta_\alpha
\end{equation}
in which $\mathbf P^\star=\{P^*_{\alpha\beta}\}_{\alpha,\beta}$
is an $\ell\times\ell$-matrix function, with entries in
$\Q[t](\l)$ (polynomial in $t$ and rational in $\l$). Their
complexity (degree and size) are bounded by
Proposition~\ref{prop:size-division}.

Since $\d(H\omega_\alpha)=H\mu_\alpha+\d H\land\omega_\alpha$
and the entries of the matrix function $\mathbf P^\star(t)$ do
not depend on $x$, the previous identities can be transformed
to the form
\begin{equation*}
    \d \biggl(H\omega_\alpha-\sum_\beta (P^\star_{\alpha\beta}\circ H)
    \cdot\omega_\beta\biggr)=-\d
    H\land(-\omega_\alpha-\eta_\alpha),
\end{equation*}
where $\d=\pd{}{x_1}\d x_1+\pd{}{x_2}\d x_2$.

Now one can choose any multiindex $s\in\Z_+^2$ and apply the
Gelfand--Leray formula with $\alpha=s$ to conclude that the
partial derivative of the matrix function $HX-\mathbf P^\star
X$ with respect to $\l_s$ is equal to the period matrix of the
forms $-x^s(\omega_\alpha+\eta_\alpha)$:
\begin{equation}\label{pf-0}
    \frac{\partial}{\partial
    \l_s}\oint_{\delta}\biggl(H\,\omega_\alpha-\sum_\beta
    (P^\star_{\alpha\beta}\circ H)\cdot\omega_\beta\biggr)=-\oint_\delta
    x^s(\omega_\alpha+\eta_\alpha).
\end{equation}

It remains to note the polynomial (in $x$) 1-forms
$x^s(\omega_\alpha+\eta_\alpha)$ can be expanded as
combinations of the basic forms $\omega_\alpha$ with
coefficients in $\Q[t](\l)$ of controlled degree and size by
\eqref{rel-coh} and Proposition~\ref{prop:size-division}.
Denote by $\mathbf
P^{s}=\bigl\{P^{s}_{\alpha\beta}\bigr\}_{\alpha,\beta}$ the
corresponding matrix functions with entries from $\Q[t](\l)$:
\begin{equation*}
    x^s(\omega_\alpha+\eta_\alpha)=\sum_\beta
    (P^{s}_{\alpha\beta}\circ H)\cdot\omega_\beta+ u_{\alpha}\d H+\d
    v_\alpha.
\end{equation*}

Substituting these identities in \eqref{pf-0} and integrating them
over the cycle $\delta\subset \{H=0\}$, on which the polynomial $H$
vanishes identically, we conclude that
\begin{equation}\label{pf-2}
    \forall s=(s_1,s_2)\in\Z_+^2,\ 0\le s_{1,2}\le n-1,
    \quad\frac{\partial}{\partial
    \l_s}\bigl(\mathbf P_0^\star X\bigr)=\mathbf P_0^{s}X.
\end{equation}
Here $\mathbf P_0^\star=\mathbf P^\star(0)$ and $\mathbf
P^{s}_0=\mathbf P^{s}(0)$ are matrices with entries in
$\Q(\l)$, obtained by setting $t=0$ in their initial
expressions. This is the Pfaffian integrable system on the
projective space $\l\in\P^m$.

From Proposition~\ref{prop:size-division} it follows that the
entries of the matrix functions $\mathbf P_0^\star,\mathbf
P_0^{s}$ are in $\Q(\l)$ and have degrees growing as $\poly(n)$
and the size (complexity) at most exponential (bounded by
$2^{\poly(n)}$). This proves Theorem~\ref{thm:pf-compl}.

\bibliographystyle{amsalpha}
\def\MR{}
\bibliography{tangent16,2003,2007}
\def\BbbR{$\mathbb R$}\def\BbbC{$\mathbb
  C$}\providecommand\cprime{$'$}\providecommand\mhy{--}\font\cyr=wncyr8\def\Bb%
bR{$\mathbb R$}\def\BbbC{$\mathbb
  C$}\providecommand\cprime{$'$}\providecommand\mhy{--}\font\cyr=wncyr9\def\Bb%
bR{$\mathbb R$}\def\BbbC{$\mathbb
  C$}\providecommand\cprime{$'$}\providecommand\mhy{--}\font\cyr=wncyr9\def\cp%
rime{$'$} \providecommand{\bysame}{\leavevmode\hbox
to3em{\hrulefill}\thinspace}
\providecommand{\MR}{\relax\ifhmode\unskip\space\fi MR }
\providecommand{\MRhref}[2]{%
  \href{http://www.ams.org/mathscinet-getitem?mr=#1}{#2}
} \providecommand{\href}[2]{#2}

\end{document}